\documentclass[leqno,11pt]{article}  
\usepackage{amsmath}
\usepackage{amssymb}

\input xy
\xyoption{all}
\usepackage[german,english]{babel}
\usepackage{amsthm}
\usepackage{xypic}  

\title{On special representations of $p$-adic reductive groups}

\author{\textsc{Elmar Grosse-Kl\"onne}}
\date{}

\theoremstyle{plain} 
\newtheorem{satz}{Theorem}[section]  
\newtheorem{kor}[satz]{Corollary}  
\newtheorem{lem}[satz]{Lemma}  
\newtheorem{pro}[satz]{Proposition}  
 
\newcommand{\ke}{\mbox{\rm Ker}}  

\newcommand{\kara}{\mbox{\rm char}}  

\theoremstyle{remark}

\theoremstyle{definition}

\begin{document}
\maketitle

\footnote[0]{MSC Classification 22E50, 11S99

Keywords: $p$-adic reductive group, modular representation, Steinberg representation}

\begin{abstract} Let $F$ be a non-Archimedean locally compact field, let $G$ be a split connected reductive group over $F$. For a parabolic subgroup $Q\subset G$ and a ring $L$ we consider the $G$-representation on the $L$-module$$(*)\quad\quad\quad\quad C^{\infty}(G/Q,L)/\sum_{Q'\supsetneq Q}C^{\infty}(G/Q',L).$$Let $I\subset G$ denote an Iwahori subgroup. We define a certain free finite rank $L$-module ${\mathfrak M}$ (depending on $Q$; if $Q$ is a Borel subgroup then $(*)$ is the Steinberg representation and ${\mathfrak M}$ is of rank one) and construct an $I$-equivariant embedding of $(*)$ into $C^{\infty}(I,{\mathfrak M})$. This allows the computation of the $I$-invariants in $(*)$. We then prove that if $L$ is a field with characteristic equal to the residue characteristic of $F$ and if $G$ is a classical group, then the $G$-representation $(*)$ is irreducible. This is the analog of a theorem of Casselman (which says the same for $L={\mathbb C}$); it had been conjectured by Vign\'{e}ras.

Herzig (for $G={\rm GL}_n(F)$) and Abe (for general $G$) have given classification
theorems for irreducible
admissible modulo $p$ representations of $G$ in terms of supersingular
representations. Some of their arguments rely on the present work.
\end{abstract}

%

\begin{center}{\bf Introduction}\end{center}

Let $F$ be a non-Archimedean locally compact field with ring of integers ${\mathcal O}_F$ and residue field $k_F$. Let $G$ be a connected split reductive group over $F$. Let $T$ be a split maximal torus, $N\subset G$ its normalizer and $W=N/T$, the corresponding Weyl group. Let $\Phi\subset X^*(T)$ be the set of roots, let $\Phi^+\subset\Phi$ be the set of positive roots with respect to a Borel subgroup $P$ containing $T$ and let $\Delta\subset\Phi^+$ be the corresponding set of simple roots. For a subset $J\subset\Delta$ let $W_J\subset W$ denote the subgroup generated by the simple reflections associated with the elements of $J$. Let $P_J$ denote the parabolic subgroup generated by $P$ and by representatives (in $N$) of the elements of $W_J$. Any parabolic subgroup of $G$ is conjugate to $P_J$ for some $J$. For a ring $L$ (commutative, with $1\in L$) we call the $G$-representation$${\rm Sp}_J(G,L)=\frac{C^{\infty}(G/P_J,L)}{\sum_{\alpha\in\Delta-J}C^{\infty}(G/P_{J\cup\{\alpha\}},L)}$$the $J$-special representation of $G$ with coefficients in $L$. For $J=\emptyset$ this is the Steinberg representation of $G$ with coefficients in $L$. By an old theorem of Casselman, the representations ${\rm Sp}_J(G,{\mathbb C})$ are irreducible for all $J$, they form the irreducible constituents, each with multiplicity one, of $C^{\infty}(G/P,{\mathbb C})$. Published proofs of this irreducibility use techniques specific for the coefficient field $L={\mathbb C}$, see \cite{bowa} ch. X, Theorem 4.11 or \cite{lau} Theorem 8.1.2. For $L$ a field of characteristic $\ell\ne p={\rm char}(k_F)$ it is known that the irreducibility of say ${\rm Sp}_{\emptyset}(G,L)$ depends on $\ell$. See e.g. \cite{virep}, Chapitre III, Th\'{e}or\`{e}me 2.8 (b).

In this paper we investigate the representation ${\rm Sp}_J(G,L)$ for arbitrary coefficient rings $L$ (and on the way obtain results previously unknown even for $L={\mathbb C}$). We need the $L$-module$${\mathfrak M}_J(L)=\frac{L[W/W_J]}{\sum_{\alpha\in\Delta-J}L[W/W_{J\cup\{\alpha\}}]}.$$ Let $I\subset G$ be an Iwahori subgroup adapted to $P$, i.e. such that we have an Iwahori decompositon $G=\bigcup_{w\in W}IwP$. Our first main theorem is the following (Theorem \ref{embedding}), which even for $L={\mathbb C}$ seems to have been unknown before:\\

{\bf Theorem A:} {\it There exists an $I$-equivariant embedding$${\rm Sp}_J(G,L)\stackrel{}{\hookrightarrow} C^{\infty}(I,{\mathfrak M}_J(L));$$its formation commutes with base changes in $L$.}\\

Using the decomposition $G/P_J=\cup_{w\in W/W_J}IwP_J/P_J$ and its analog for the $P_{J\cup\{\alpha\}}$, the proof of Theorem 1 is reduced to the proof of exactness of a certain natural sequence\begin{gather}\bigoplus_{{\alpha}\in\Delta-J\atop w\in W/W_{J\cup\{{\alpha}\}}}C^{\infty}(I/I\cap wP_{J\cup\{{\alpha}\}}w^{-1},L)\stackrel{}{\longrightarrow}\bigoplus_{w\in W/W_{J}}C^{\infty}(I/I\cap wP_{J}w^{-1},L)\stackrel{}{\longrightarrow}C^{\infty}(I,{\mathfrak M}_J(L))\label{introsec}\end{gather}(Proposition \ref{keyiwa}). This exactness proof proceeds by induction along a certain filtration of (\ref{introsec}). The key to defining this filtration is to consider certain subsets of $\Phi$ which we call $J$-quasi-parabolic: a subset $D\subset\Phi$ is called $J$-quasi-parabolic if $\prod_{\alpha\in D}U_{\alpha}$ is the intersection of unipotent radicals of parabolic subgroups which are $W$-conjugate to $P_J$. Here $U_{\alpha}\subset G$ denotes the root subgroup associated to $\alpha$. For such $D$ we define a subset $W^J(D)$ of $W/W_J$ as consisting of those classes $wW_J$ for which $\prod_{\alpha\in D}U_{\alpha}$ is contained in the unipotent radical of the parabolic subgroup opposite to $wP_Jw^{-1}$. Fixing a size-increasing enumeration of all $J$-quasi-parabolic subsets $D$, the corresponding $W^J(D)$'s give the said filtration of (\ref{introsec}). The exactness of (\ref{introsec}) is then reduced to the exactness, for any $D$, of$$\bigoplus_{{\alpha}\in\Delta-J}L[W^{J\cup\{{\alpha}\}}(D)]\stackrel{}{\longrightarrow}L[W^J(D)]\stackrel{}{\longrightarrow}{\mathfrak M}_J(L)$$(Proposition \ref{parml}), a purely combinatorial fact on finite reflection groups. We mention that if $L$ is a complete field extension of $F$, Theorem 1 holds verbatim, with the same proof, for the corresponding representations on spaces of locally analytic (rather than locally constant) functions.

A vigorously emerging subject in current $p$-adic number theory is the smooth representation theory of $p$-adic reductive groups, like $G$, on $\overline{\mathbb F}_p$-vector spaces. So far, the research has focused mostly on the case $G={\rm GL}_2(F)$, for finite extensions $F$ of ${\mathbb Q}_p$, but even for those $G$ the theory turns out to be fairly complicated and is far from being well understood. However, it already becomes quite clear that a good understanding of the theory depends crucially on a good understanding of the functor {\it taking invariants under a (pro-$p$-)Iwahori-subgroup}. At present there is literally no general technique available to compute this functor. For example, although  Vign\'{e}ras had proved the irreducibility of the Steinberg representation of our $G$'s in characteristic $p$, the space of its (pro-$p$-)Iwahori invariants was not known (except for $G={\rm GL}_2(F)$); this was the motivating problem for our investigations. 

As an immediate consequence of Theorem 1 we obtain that the submodule of
$I$-invariants ${\rm Sp}_J(G,L)^I$ is free of rank at most the rank of
${\mathfrak M}_J(L)$, i.e. ${\rm rk}_L({\rm Sp}_J(G,L)^I)\le{\rm
  rk}_L({\mathfrak M}_J(L))$, as was conjectured by Vign\'{e}ras
\cite{vig}. The reverse inequality ${\rm rk}_L({\rm Sp}_J(G,L)^I)\ge{\rm
  rk}_L({\mathfrak M}_J(L))$ follows easily by summing over all $J$, using
that $\sum_J{\rm rk}_L({\mathfrak M}_J(L))=|W|$. Thus, ${\rm Sp}_J(G,L)^I$ is
free of rank equal to the rank of ${\mathfrak M}_J(L)$, for any $L$. (For
example, we obtain that the module of $I$-invariants in the Steinberg
representation is free of rank one.) In particular, using Lemma 6.18 of \cite{pas}:\\

{\bf Corollary B:} {\it The $G$-representation ${\rm Sp}_J(G,L)$ is
  admissible, for any $J$ and any $L$.}\\

(Corollary 2 also follows from Proposition 2.2.13 of
\cite{eme} and the admissibility of ${C}^{\infty}(G/P_J,L)$.) The reductive
group underlying $G$ can be defined over ${\mathcal O}_F$; as such we denote
it by ${\mathcal G}_{x_0}$. Its group ${\mathcal G}_{x_0}({\mathcal O}_F)$ of
${\mathcal O}_F$-rational points is a subgroup of $G$, let
$\overline{G}={\mathcal G}_{x_0}(k_F)$ denote the group of $k_F$-rational
points of ${\mathcal G}_{x_0}$. Its root system is the same as that of $G$. We may copy the definition of the $G$-representations ${\rm Sp}_J(G,L)$ to define $\overline{G}$-representations ${\rm Sp}_J(\overline{G},L)$, for all $J\subset\Delta$ (replace locally constant functions on $G$ by functions on $\overline{G}$). Let $\overline{P}\subset\overline{G}$ denote the Borel subgroup obtained by reduction of $I\subset {\mathcal G}_{x_0}({\mathcal O}_F)$. Then using Theorem 1 we find a canonical identification (Proposition \ref{ux0in}):\begin{gather}{\rm Sp}_J({G},L)^{I}={\rm Sp}_J(\overline{G},L)^{\overline{P}}.\label{introux0}\end{gather}

Our second main theorem is concerned with the case where $L$ is a field with
$p=\kara(L)=\kara(k_F)$. We ask whether ${\rm Sp}_J(G,L)^I$ is irreducible as
a module under the Iwahori Hecke algebra ${\mathcal H}(G,I)$. We may view
${\rm Sp}_J(G,L)^I={\rm Sp}_J(\overline{G},L)^{\overline{P}}$ as a module
under the Hecke algebra ${\mathcal H}(\overline{G},\overline{P})$. In a first
step we show (Proposition \ref{indeco}) that each ${\mathcal
  H}(\overline{G},\overline{P})$-submodule of ${\rm Sp}_J(G,L)^I={\rm
  Sp}_J(\overline{G},L)^{\overline{P}}$ contains the class of the
characteristic function $\chi_{Iw_{\Delta}P_J}$ of the subset
$Iw_{\Delta}P_J\subset G$; here $w_{\Delta}\in W$ denotes the longest
element. This follows from explicit formulae for the action on ${\rm
  Sp}_J(\overline{G},L)^{\overline{P}}$ of the Hecke operators associated to
simple reflections (these formulae boil down to the Bruhat decomposition of
$\overline{G}$ and require our assumption $p=\kara(L)=\kara(k_F)$), together
with a combinatorial lemma (Lemma \ref{weyllem1}) on $W$. In a second step we
need to show that the class of $\chi_{Iw_{\Delta}P_J}$ generates ${\rm
  Sp}_J(G,L)^I$ as an ${\mathcal H}(G,I)$-module. We can prove this if $\Phi$
contains no exceptional factor, i.e. if all the irreducible
factors of the root system $\Phi$ belong to the infinite series $A$, $B$, $C$
or $D$. Our argument uses a combinatorial result, Proposition \ref{weyllem2},
on the weak (left)ordering of $W$ (an ordering weaker than the Bruhat
ordering) which we can prove only for such root systems. Proposition
\ref{weyllem2} may also hold true for the root systems of type $E_6$ or $E_7$
(if so we would get the irreducibility result in these cases too), but certainly fails for the root systems of the types $E_8$, $F_4$ and $G_2$. Thus, in these cases another argument (for the generation of ${\rm Sp}_J(G,L)^I$ by $\chi_{Iw_{\Delta}P_J}$) would be needed. In conclusion, what we prove is (Theorem \ref{hecirr}):\\

{\bf Theorem C:} {\it If $L$ is a field with $\kara(L)=\kara(k_F)$ and if the
  root system $\Phi$ contains no exceptional factor then the ${\mathcal H}(G,I;L)$-module ${\rm Sp}_J({G},L)^I$ is irreducible.}\\

Let $I_1\subset I$ denote the pro-$p$-Iwahori subgroup inside $I$. The $G$-representation ${\rm Sp}_J({G},L)$ is generated by ${\rm Sp}_J(G,L)^I={\rm Sp}_J(G,L)^{I_1}$ (see \cite{vig}). As any smooth representation of a pro-$p$-group on a non-zero vector space in characteristic $p$ admits a non-zero invariant vector, we obtain, as a corollary of Theorem C, the analog of Casselman's theorem for a field $L$ with $p=\kara(L)=\kara(k_F)$ if $G$ is a classical group (of course, this analog implies and gives a new, purely algebraic proof of Casselman's theorem) (Corollary \ref{girr}, Corollary \ref{johoe}):\\

{\bf Theorem D:} {\it If $L$ is a field with $\kara(L)=\kara(k_F)$ and if $\Phi$ contains no exceptional factor then the $G$-representation ${\rm Sp}_J({G},L)$ is irreducible. The ${\rm Sp}_J({G},L)$ for the various $J$ form the irreducible constituents, each one occuring with multiplicity one, of $C^{\infty}(G/P,L)$.}\\

Theorem 4 had been conjectured by Vign\'{e}ras (see \cite{vig} section 5, Remarque 2) (without the
restriction on $\Phi$), and, as indicated above, she had proven the
irreducibility of the Steinberg representation ${\rm
  Sp}_{\emptyset}({G},L)$. After we had obtained Theorem 4 it had been
generalized by Florian Herzig \cite{her} to general (split reductive) groups
$G$ over a finite extension $F$ of ${\mathbb Q}_p$. Like ours, his proof
relies on the identification (\ref{introux0}) and on Proposition \ref{indeco} below, but then it follows another strategy; in particular, it does not reprove or generalize Theorem 3.

Assuming the results of the present paper, Florian Herzig \cite{her} (for $G={\rm GL}_n(F)$) and Noriyuki Abe \cite{abe} (who
generalized Herzig's method to general split $G$) classify irreducible
admissible representations of $G$ over $L$ in terms of supersingular
representations; here $G$ is a split connected reductive group $G$ over a
finite extension $F$ of ${\mathbb Q}_p$ and $L$ is an algebraically closed
field $L$ with $\kara(L)=\kara(k_F)=p$. More specifically, our results
(e.g. Corollary 2, formula (\ref{introux0}), Proposition \ref{indeco}) are
indispensable for proving e.g. the irreducibility of the representations considered in these papers.

It is a great pleasure to express my deep gratitude to Marie-France
Vign\'{e}ras. She suggested the problem of computing the Iwahori invariants in
$p$-modular Steinberg representations: this was the origin of the present
work. Later she gave helpful comments on a preliminary version of this
paper. I am extremely grateful to Peter Schneider. Having explained to him an
unnecessarily complicated proof of Theorem 1, valid only in a restricted
setting, he insisted on getting a better conceptual understanding. His
numerous suggestions were decisive for approaching Theorem 1 in the correct
context and for discovering the proof in its full generality. He also outlined
some possible further developments. I thank Florian Herzig for his very
careful reading of the manuscript and the numerous email exchanges which we
had about it. The referees wrote detailed and helpful reports for which I
am very grateful. I thank the Deutsche Forschungs Gemeinschaft (DFG) as part of this work was done while I was supported by the DFG as a Heisenberg fellow.

\tableofcontents

\section{Reflection groups}
\label{weysec}

In this section we collect some results on finite reflection groups. Proposition \ref{parml} will be needed for Theorem \ref{embedding}, the embedding of ${\rm Sp}_J({G},L)$ into $C^{\infty}(I,{\mathfrak M}_J(L))$. Lemma \ref{weyllem1} will be needed for Proposition \ref{indeco} which concerns the ${\mathcal H}(\overline{G},\overline{P};L)$-module structure of ${\rm Sp}_J({G},L)^{I}$, and Corollary \ref{weyllem3} will be needed for the proof of Theorem \ref{hecirr} on the irreducibility of ${\rm Sp}_J({G},L)^I$ as a ${\mathcal H}(G,I;L)$-module.

Consider a reduced root system $\Phi$ and let $W$ be its corresponding Weyl
group. Fix a system $\Delta\subset\Phi$ of simple roots and denote by
$\Phi^+\subset \Phi$ the corresponding set of positive roots. Let
$\Phi^-=\Phi-\Phi^+$.  For $\alpha\in\Phi$ let $s_{\alpha}\in W$ denote the
associated reflection. Let $\ell(.):W\to{\mathbb Z}_{\ge0}$ be the length
function with respect to $\Delta$. For a subset $J\subset\Delta$ let
$W_J\subset W$ be the subgroup generated by all $s_{\alpha}$ for $\alpha\in
J$. We denote by $w_{\Delta}\in W$ resp. $w_J\in W_J$ the respective longest
elements. Let$$\Phi_J(1)=\Phi^--(\Phi^-\cap W_J.J)$$where
$W_J.J=\{w\alpha\,|\,w\in W_J, \alpha\in J\}\subset\Phi$ is the sub-root
system generated by $J$. For $w\in W$ we then define the
subset$$\Phi_J(w)=w\Phi_J(1)$$of $\Phi$. It depends only on the class of $w$
in $W/W_J$. Observe $\Phi_{J'}(w)\subset\Phi_{J}(w)$ for $J\subset J'$. We say
that a subset $D\subset\Phi$ is $J${\it -quasi-parabolic} if it is the
intersection of subsets $\Phi_J(w)$ for some (at least one) $w\in
W$. Let$$W^J=\{w\in W\,\,|\,\,w(J)\subset\Phi^+\}.$$It is well known
(cf. e.g. \cite{hum} Proposition 1.10 (c)) that this is a set of representatives for $W/W_J$ and can alternatively be described as\begin{gather}W^J=\{w\in W\,\,|\,\,\ell(ws_{\alpha})>\ell(w)\mbox{ for all }{\alpha}\in J\}.\label{lenchar}\end{gather}
For a subset $D\subset \Phi$ let$$W^J(D)=\{w\in W^J\,\,|\,\,D\subset\Phi_J(w)\}.$$Let$$V^J=W^J-\bigcup_{{\alpha}\in\Delta-J}W^{J\cup\{{\alpha}\}}.$$Then $W=\cup_{J\subset\Phi}V^J$ (disjoint union). We have$$V^J=\{w\in W^J\,\,|\,\,w(\Delta-J)\subset\Phi^-\}.$$

\begin{lem}\label{hilfe} For $J\subset J'$ and $w\in W^{J'}$ we have $\Phi_J(w)-\Phi_{J'}(w)\subset\Phi^-$.
\end{lem}

{\sc Proof:} Each element in $\Phi_J(w)-\Phi_{J'}(w)=w(\Phi_J(1)-\Phi_{J'}(1))$ can be written as $w(\sum_{\nu}-\alpha_{\nu})$ with certain $\alpha_{\nu}\in J'$. As $w\in W^{J'}$ the claim follows.\hfill$\Box$\\ 

For the proof of Proposition \ref{parml} below and then for later use it is convenient to make the following definition:\\

{\bf Definition:} For $w\in W$ let $(w)^J$ denote the unique element of $W^J$ with $(w)^JW_J=wW_J$. Thus, $(.)^J$ is the projection from $W$ onto the first factor in the direct product decomposition $W=W^JW_J$. Loosely speaking, applying $(.)^J$ means cutting off $W_J$-factors on the right hand side. 

\begin{lem}\label{warmup} (a) For any $w\in W$ we have $\ell(w)\ge \ell((w)^J)$.\\(b) For $w_1\in W^J$ and $w_2\in W_J$ we have $\ell(w_1w_2)=\ell(w_1)+\ell(w_2)$.\\(c) For any $w\in W$ we have $\ell(w_{\Delta}w)=\ell(ww_{\Delta})=\ell(w_{\Delta})-l(w)$.
\end{lem}

{\sc Proof:} Any $v\in W^J$ is the unique element of minimal length in the set of representatives for the coset $vW_J$; this gives (a). For the easy statements (b) and (c) see \cite{hum} Theorem 1.8 and Proposition 1.10.\hfill$\Box$\\

Let $L$ be a ring. For a set $S$ let $L[S]$ denote the free $L$-module with basis $S$.\\

{\bf Definition:} We define the $L$-module ${\mathfrak M}_J(L)$ and the $L$-linear map $\nabla$ by the exact sequence of $L$-modules\begin{gather}\bigoplus_{{\alpha}\in\Delta-J}L[W^{J\cup\{{\alpha}\}}]\stackrel{\partial}{\longrightarrow}L[W^J]\stackrel{\nabla}{\longrightarrow}{\mathfrak M}_J(L){\longrightarrow}0\label{defmjl}\end{gather}where for $w\in W^{J\cup\{{\alpha}\}}$ we set $$\partial(w)=\sum_{w'\in W^J\atop w'W_J\subset wW_{J\cup\{{\alpha}\}}}w'.$$

\begin{pro}\label{parml} (a) ${\mathfrak M}_J(L)$ is $L$-free of rank $|V^J|$,
  and $\nabla$ induces a bijection between ${V}^J$ and an $L$-basis of
  ${\mathfrak M}_J(L)$. We have ${\mathfrak M}_J(L')={\mathfrak M}_J(L)\otimes_LL'$ for any ring morphism $L\to L'$.

(b) Let $D\subset\Phi$ be a $J$-quasi-parabolic subset. We have $\partial(\oplus_{{\alpha}\in\Delta-J}L[W^{J\cup\{{\alpha}\}}(D)])\subset L[W^J(D)]$, and the sequence$$\bigoplus_{{\alpha}\in\Delta-J}L[W^{J\cup\{{\alpha}\}}(D)]\stackrel{\partial^D}{\longrightarrow}L[W^J(D)]\stackrel{\nabla^D}{\longrightarrow}{\mathfrak M}_J(L)$$obtained by restricting (\ref{defmjl}) is exact.
\end{pro}

{\sc Proof:} For $w\in W^{J\cup\{{\alpha}\}}$ and $w'\in W^J$ with $w'W_J\subset wW_{J\cup\{\alpha\}}$ we have $\Phi_{J\cup\{\alpha\}}(w)=\Phi_{J\cup\{\alpha\}}(w')\subset\Phi_J(w')$. This shows $$\partial(\oplus_{{\alpha}\in\Delta-J}L[W^{J\cup\{{\alpha}\}}(D)])\subset L[W^J(D)],$$for any subset $D$ of $\Phi$.

{\it First Step:} Let $D\subset\Phi^+$ be a subset. Define ${\mathfrak M}_{J,D}(L)$ and $\tilde{\nabla}^D$ by the exact sequence$$\bigoplus_{{\alpha}\in\Delta-J}L[W^{J\cup\{{\alpha}\}}(D)]\stackrel{\partial^D}{\longrightarrow}L[W^J(D)]\stackrel{\tilde{\nabla}^D}{\longrightarrow}{\mathfrak M}_{J,D}(L){\longrightarrow}0.$$Let $V^J(D)={V}^J\cap W^J(D)$.

{\it Claim:} For all $\ell$ and all $w\in W^J(D)$ with $\ell(w)\ge\ell$ we have $\tilde{\nabla}^D(w)\in \tilde{\nabla}^D(L[V^J(D)])$.

We prove this by descending induction on $\ell$. Suppose we are given such a $w\in W^J(D)$ with $\ell(w)\ge\ell$. If $w\in {V}^J$ we are done. Otherwise there is some ${\alpha}\in\Delta-J$ with $w\in W^{J\cup\{{\alpha}\}}$. By Lemma \ref{hilfe} we have $\Phi_J(w)-\Phi_{J\cup\{\alpha\}}(w)\subset\Phi^-$, thus our assumption $D\subset\Phi^+$ implies even $w\in W^{J\cup\{{\alpha}\}}(D)$. For all $w'\in W^J-\{w\}$ with $w'W_J\subset wW_{J\cup\{{\alpha}\}}$ we have $\ell(w')>\ell(w)$ (because $w'W_J\subset wW_{J\cup\{{\alpha}\}}$ implies $w'W_{J\cup\{{\alpha}\}}=wW_{J\cup\{{\alpha}\}}$, but in view of (\ref{lenchar}) we know that $w$ is the unique element of $wW_{J\cup\{{\alpha}\}}$ of minimal length). Moreover we have $w'\in W^J(D)$ (as noted at the beginning of this proof), thus by induction hypothesis we get $\tilde{\nabla}^D(w')\in \tilde{\nabla}^D(L[V^J(D)])$ for all such $w'$. Now $$w=\partial^D(w)-\sum_{w'\in W^J-\{w\}\atop w'W_J\subset wW_{J\cup\{{\alpha}\}}}w'$$(inside $L[W^J(D)]$) which shows $\tilde{\nabla}^D(w)\in \tilde{\nabla}^D(L[V^J(D)])$, as desired.

The claim is proved. In particular, setting $\ell=0$, we get $\tilde{\nabla}^D(L[V^J(D)])={\mathfrak M}_{J,D}(L)$.

{\it Second Step:} Here we prove (a). That the image of ${V}^J$ generates the
$L$-module ${\mathfrak M}_J(L)$ follows from the first step (with
$D=\emptyset$ there). The base change property ${\mathfrak
  M}_J(L')={\mathfrak M}_J(L)\otimes_LL'$ follows from the definition of ${\mathfrak
  M}_J(.)$ and from right exactness of taking tensor products. To see that
the image of ${V}^J$ in ${\mathfrak M}_J(L)$ remains linearly independent we
first consider the case $L={\mathbb Q}$; then our task is to show $\dim_{\mathbb Q}{\mathfrak M}_J(\mathbb Q)\ge|V^J|$.

By definition, the ${\mathbb Q}$-vector spaces ${\mathbb Q}[W^J]$ and ${\mathbb Q}[W^{J\cup\{\alpha\}}]$ come with the distinguished bases $W^J$ and $W^{J\cup\{\alpha\}}$, hence with isomorphisms with their duals ${\mathbb Q}[W^J]\cong {\mathbb Q}[W^J]^*$ and ${\mathbb Q}[W^{J\cup\{\alpha\}}]\cong{\mathbb Q}[W^{J\cup\{\alpha\}}]^*$. One easily checks that under these identifications, the map$${\mathbb Q}[W^J]\stackrel{\partial^*}{\longrightarrow}\bigoplus_{{\alpha}\in\Delta-J}{\mathbb Q}[W^{J\cup\{{\alpha}\}}]$$dual to $\partial$ is given as follows: for $w'\in W^J$ the $\alpha$-component of $\partial^*(w')$ is the unique $w\in W^{J\cup\{{\alpha}\}}$ with $w'W_{J\cup\{\alpha\}}= wW_{J\cup\{{\alpha}\}}$. For $w'\in V^J$ put$$\sigma(w')=\sum_{v\in W_{\Delta-J}}(-1)^{\ell(v)}(w'v)^J\in {\mathbb Q}[W^J].$$The definition of $V^J$ shows that for each $w'\in V^J$ and each $v\in W_{\Delta-J}$ different from the neutral element we have $\ell(w')>\ell(w'v)\ge\ell((w'v)^J)$. Therefore the set $$\{\sigma(w')\,|\,w'\in V^J\mbox{ and }\ell(w')=\ell\}$$remains linearly independent in $$\frac{{\mathbb Q}[W^J]}{{\mathbb Q}[\{w'\in W^J\,|\,\ell(w')<\ell\}]}$$ for any $\ell\in{\mathbb N}$. An induction then shows that the set $\{\sigma(w')\,|\,w'\in V^J\}$ is linearly independent in ${\mathbb Q}[W^{J}]$ (under the projection ${\mathbb Q}[W^{J}]\to {\mathbb Q}[V^{J}]$ it even maps bijectively onto a ${\mathbb Q}$-basis of ${\mathbb Q}[V^{J}]$). On the other hand, for any $\alpha\in \Delta-J$ we have $W_{\Delta-J}=(W_{\Delta-J})^{\alpha}\coprod(W_{\Delta-J})^{\alpha}s_{\alpha}$ (we extrapolate to $W_{\Delta-J}$ the definitions given for $W$, i.e. $(W_{\Delta-J})^{\alpha}$ is the set of canonical representatives for $W_{\Delta-J}/W_{\{\alpha\}}$). Therefore the above description of ${\partial^*}$ shows that $\sigma(w')\in {\rm ker}(\partial^*)$ for all $w'\in V^J$. We obtain $\dim_{\mathbb Q}{\mathfrak M}_J({\mathbb Q})={\rm dim}_{\mathbb Q}{\rm coker}(\partial)={\rm dim}_{\mathbb Q}{\rm ker}(\partial^*)\ge|V^J|$, as desired.

We have proven that the image of ${V}^J$ in ${\mathfrak M}_J({\mathbb Q})$ is
a ${\mathbb Q}$-basis of ${\mathfrak M}_J({\mathbb Q})$. Since the image of ${V}^J$
in ${\mathfrak M}_J({\mathbb Z})$ generates ${\mathfrak M}_J({\mathbb Z})$ as
an abelian group, and as ${\mathfrak M}_J({\mathbb Q})= {\mathfrak
  M}_J({\mathbb Z})\otimes{\mathbb Q}$, it follows that ${\mathfrak
  M}_J({\mathbb Z})$ is torsion free and that the image of ${V}^J$ in
${\mathfrak M}_J({\mathbb Z})$ is a ${\mathbb Z}$-basis. By the base change
property it follows that ${\mathfrak M}_J({\mathbb Q})$ is $L$-free for any
$L$, with the image of ${V}^J$ as an $L$-basis.

{\it Third Step:} Here we prove (b). As $D$ is $J$-quasi-parabolic we find some $w\in W$ with $wD\subset\Phi^+$. We have a commutative diagram$$\xymatrix{\bigoplus_{\alpha\in \Delta-J}L[W^{J\cup\{\alpha\}}(D)]\ar[d]^{\cong}\ar[r]^{\quad\quad\quad{\partial}^D}& L[W^{J}(D)]\ar[d]^{\cong}\ar[r]^{\quad{\nabla}^D}&{\mathfrak M}_J(L)\ar[d]^{\cong}\\{\bigoplus}_{{\alpha}\in\Delta-J}L[W^{J\cup\{\alpha\}}(wD)]\ar[r]^{\quad\quad\quad\quad{\partial}^{wD}}&L[W^{J}(wD)]\ar[r]^{\quad{\nabla}^{wD}}&{\mathfrak M}_J(L)}$$where the second and the third (resp. the first) vertical isomorphism is induced by the bijection $W^J\to W^J$, $w'\mapsto (ww')^J$ (resp. $W^{J\cup\{\alpha\}}\to W^{J\cup\{\alpha\}}$, $w'\mapsto (ww')^{J\cup\{\alpha\}}$). Therefore we may assume from the beginning that $D\subset\Phi^+$. It suffices to see that the natural map ${\mathfrak M}_{J,D}(L)\to{\mathfrak M}_J(L)$ is injective. By (a) we know that the image of ${V}^J$, hence in particular the image of $V^J(D)$ in ${\mathfrak M}_J(L)$ is linearly independent. Together with the result of the first step this shows the wanted injectivity of ${\mathfrak M}_{J,D}(L)\to{\mathfrak M}_J(L)$.\hfill$\Box$\\

{\bf Definition:} We write $S=\{s_{\alpha}\,|\,\alpha\in\Delta\}$. Consider the following partial ordering $<_J$ on $W^J$. For $w, w'\in W^J$ we write $w<_Jw'$ if there are $s_1,\ldots,s_r\in S$ such that, setting $w^{(i)}=(s_i\cdots s_1w)^J$ for $0\le i\le r$, we have $\ell(w^{(i-1)})<\ell(w^{(i)})$ for all $i\ge 1$, and $w^{(r)}=w'$.

\begin{lem}\label{weylem} Let $w\in W^J$ and $s\in S$.\\(a) If $w<_J(sw)^J$
  then we have $\ell(w)<\ell(sw)$.\\(b) $\ell(w)<\ell(sw)$ and $w\ne(sw)^J$
  together imply $sw\in W^J$, hence $w<_J(sw)^J=sw$. We have $$(sw)^J=w\quad\mbox{ or }\quad(sw)^J=sw.$$\\(c) $$(sw)^J<_Jw\quad\quad \Leftrightarrow\quad\quad  \ell((sw)^J)<\ell(w) \quad \quad \Leftrightarrow\quad\quad  \ell(sw)<\ell(w).$$(d) Let $u\in W$. If $w_Jw_{\Delta}<_{\emptyset}uw_{\Delta}$ then $u\in W_J$.\\(e) There exists a unique maximal element $z^J\in W^J$ for the ordering $<_J$; it lies in $V^J$. We have $z^J=w_{\Delta}w_J$. For any $u\in W$ such that $z^J\le_{\emptyset}u$ and for any $s\in S$ with $\ell(sz^J)<\ell(z^J)$ we have $\ell(su)<\ell(u)$.\\(f) If $w\in V^J$ and $\ell((sw)^J)>\ell(w)$ then $(sw)^J\in V^J$.\end{lem}

{\sc Proof:} (a) We have $\ell(w)<\ell((sw)^J)\le \ell(sw)$ where the first inequality follows from the definition of $<_J$ and the second one from Lemma \ref{warmup} (a) (applied to $sw$).

To prove (b) assume $\ell(w)<\ell(sw)$ and $sw\notin W^J$. Then we find some ${\alpha}\in J$ with $\ell(sws_{\alpha})=\ell(sw)-1=\ell(w)$. Take a reduced expression $w=\sigma_1\cdots\sigma_r$ with $\sigma_i\in S$. By the deletion condition for Weyl groups we get a reduced expression for $sws_{\alpha}$ by deleting some factors in the string $s\sigma_1\ldots\sigma_rs_{\alpha}$. Namely, as $\ell(sws_{\alpha})=\ell(w)$, exactly two factors must be deleted. If $s$ remained this would mean $\ell(ws_{\alpha})<\ell(w)$, contradicting $w\in W^J$. If $s_{\alpha}$ remained this would mean $\ell(sw)<\ell(w)$, contradicting our hypothesis. Thus $sws_{\alpha}=w$, i.e. $w=(sw)^J$.

(c) First assume $\ell(sw)<\ell(w)$. Then we get $\ell((sw)^J)<\ell(w)$ from Lemma \ref{warmup} (a) (applied to $sw$). As $(s(sw)^J)^J=w^J=w$ we get $(sw)^J<_Jw$ from the definition of $<_J$. If on the other hand we have $\ell(sw)>\ell(w)$ then we cannot have  $(sw)^J<_Jw$ at the same time, as follows from (b). We have shown the equivalence of the outer statements. Since by (b) we always have $(sw)^J=w$ or $(sw)^J=sw$ they are equivalent with the middle statement.

(d) Letting $v=uw_J$, the statement $u\in W_J$ is equivalent with the statement $v\in W_J$. Consider the following
chain of equalities$$\ell(w_{\Delta})=\ell(vw_Jw_{\Delta})+\ell(w_J
v^{-1})=\ell(v)+\ell(w_Jw_{\Delta})+\ell(w_Jv^{-1})=\ell(w_Jw_{\Delta})+\ell(w_J)=\ell(w_{\Delta}).$$Here
the second equality follows from our hypothesis
$w_Jw_{\Delta}<_{\emptyset}uw_{\Delta}=vw_Jw_{\Delta}$. The third equality
follows from the conjunction of all the other equalities (and the equality of
the extreme terms in the chain). But this third equality says
$\ell(v)+\ell(w_Jv^{-1})=\ell(w_J)$ which implies $v\in W_J$, because no
reduced expression for $w_J$ contains an $s_{\alpha}$ with $\alpha\in\Delta-J$
(if it did, then, by the subword property in Coxeter groups, $s_{\alpha}$
would occur in {\it any} reduced expression of $w_J$, which is nonsense). 

As a referee pointed out, statement (d) follows alternatively from well
known results on the Bruhat order, because $w_Jw_{\Delta}<_{\emptyset}uw_{\Delta}$
implies that $w_J$ is larger than $u$ in the Bruhat order.
 
(e) From Lemma \ref{warmup} (c) it follows that
$(w_{\Delta})^J=w_{\Delta}w_J$. We claim that
$z^J=(w_{\Delta})^J=w_{\Delta}w_J$ is maximal in $W^J$ with respect to $<_J$,
and is uniquely determined by this property. To see this we need to show, by
(b), that for any $w\in W^J-\{z^J\}$ there is some $s\in S$ with
$\ell(sw)>\ell(w)$ and $w\ne(sw)^J$. As $w\ne z^J=w_{\Delta}w_J$ we find $s\in
S$ with $\ell(sww_J)=\ell(ww_J)+1$, hence$$\ell(sw)\ge
\ell(sww_J)-\ell(w_J)=\ell(ww_J)+1-\ell(w_J)>\ell(w)$$where we used
$\ell(ww_J)=\ell(w)+\ell(w_J)$ as recorded in Lemma \ref{warmup} (b). If we
had $w=(sw)^J$ this would mean $sw=wu$ for some $u\in W_J$, hence
$\ell(sww_J)=\ell(wuw_J)\le \ell(ww_J)$ by Lemma \ref{warmup} (b):
contradiction ! The claim is proved.

For $\alpha\in\Delta-J$ we have $\ell(s_{\alpha}w_J)>\ell(w_J)$. Since
$w_{\Delta}=z^Jw_J=(z^{J}s_{\alpha})(s_{\alpha}w_J)$ we thus get
$\ell(z^Js_{\alpha})=\ell(w_{\Delta})-\ell(s_{\alpha}w_{J})<\ell(w_{\Delta})-\ell(w_{J})=\ell(z^J)$,
hence $z^J\in V^J$.

Finally, we have $z^J=w_{\Delta}w_J=w_{\check{J}}w_{\Delta}$
for $$\check{J}=\{\beta\in\Delta\,|\,s_{\beta}=w_{\Delta}s_{\alpha}w_{\Delta}\mbox{
  for some }\alpha \in J\}.$$Equivalently, $\check{J}=-w_{\Delta}(J)$. For
$u\in W$ such that
$z^J=w_{\check{J}}w_{\Delta}<_{\emptyset}u=(uw_{\Delta})w_{\Delta}$ we get
$uw_{\Delta}\in W_{\check{J}}$ using (d). The same argument which showed
$z^J\in V^J$ also shows that $\ell(sz^J)<\ell(z^J)$ for $s\in S$ can only happen if $s=s_{\alpha}$ for some $\alpha\in \Delta-\check{J}$. Therefore $\ell(suw_{\Delta})>\ell(uw_{\Delta})$ since $uw_{\Delta}\in W_{\check{J}}$. By Lemma \ref{warmup} (c) this means $\ell(su)<\ell(u)$.

(f) Follows from (the proof of) (c).\hfill$\Box$\\

\begin{lem}\label{weyllem1} For each $w\in V^J-\{z^J\}$ there is some $w'\in V^J$ and some $s\in S$ with $w<_Jw'$, with $\ell((sw)^J)<\ell(w)$ and with $\ell((sw')^J)\ge\ell(w')$.
\end{lem}

{\sc Proof:} Consider the set$$J'=\{\alpha\in
\Delta\,\,|\,\,\ell(s_{\alpha}w)>\ell(w)\}.$$For any given $\alpha\in\Delta$
we have $\alpha\notin J'$ if and only if $\ell((s_{\alpha}w)^J)<\ell(w)$, by
Lemma \ref{weylem}(c).

{\it Case (i): $z^Jw^{-1}\notin W_{J'}$.} As $z_J$ is maximal for the ordering
$<_J$ on $W_J$ (Lemma
\ref{weylem}(e)), we find $\sigma_1,\ldots,\sigma_r$ in $S$ with
$w<_J(\sigma_1w)^J<_J\ldots<_J(\sigma_r\cdots\sigma_1 w)^J=z^J$. Lemma
\ref{weylem}(b), applied first to $w<_J(\sigma_1w)^J$, then to
$(\sigma_1w)^J<_J(\sigma_2\sigma_1w)^J$, then to
$(\sigma_2\sigma_1w)^J<_J(\sigma_3\sigma_2\sigma_1w)^J$ etc. shows
successively that
$(\sigma_j\cdots\sigma_1 w)^J=\sigma_j\cdots\sigma_1 w$ for all $j$. We get $\sigma_r\cdots\sigma_1 w=z^J$ and $\ell(z^J)=r+\ell(w)$. Let $1\le i\le r$
be maximal such that $\sigma_i=s_{\alpha}$ for some $\alpha\in\Delta-J'$ (such
an $i$ exists since $z^Jw^{-1}\notin W_{J'}$). By Lemma \ref{weylem}(b) we
then see $w'\in W^J$ for $w'= \sigma_{i+1}\cdots\sigma_rw$. But then we
necessarily even have $w'\in V^J$. Indeed, otherwise we would have $w'\in
W^{J\cup\{\alpha\}}$ for some $\alpha\in\Delta-J$, hence
$\ell(\sigma_{i+1}\cdots\sigma_rws_{\alpha})=\ell(w's_{\alpha})>\ell(w')=\ell(w)+r-i$. On the other hand, as $w\in
V^J$ we have $\ell(w)>\ell(ws_{\alpha})$, and together we would obtain a contradiction. Thus, this $w'$ together with $s=s_{\alpha}$ is fine. 

{\it Case (ii): $z^Jw^{-1}\in W_{J'}$.} Note that this implies
$z^J\le_{\emptyset}w_{J'}w$ (because of $\ell(w_{J'}w)=\ell(w_{J'})+\ell(w)$
as follows from the definition of $J'$). Here we claim that $w'=z^J$ satisfies
the wanted conclusion. Assume on the contrary that
$\ell(s_{\alpha}z^J)<\ell(z^J)$ for all $\alpha\in \Delta-J'$. Then we also
have $\ell(s_{\alpha}w_{J'}w)<\ell(w_{J'}w)$ for all $\alpha\in
\Delta-J'$. This follows from Lemma \ref{weylem}(e) since
$z^J\le_{\emptyset}w_{J'}w$ as noted above. On the other hand
$\ell(s_{\alpha}w_{J'}w)<\ell(w_{J'}w)$ for all $\alpha\in J'$, too (again
because of $\ell(w_{J'}w)=\ell(w_{J'})+\ell(w)$), hence for all
$\alpha\in\Delta$. This means $w_{J'}w=w_{\Delta}$. But then
$w=w_{\Delta}w_{\check{J}}$ for some $\check{J}\subset\Delta$ (as in the proof
of Lemma \ref{weylem}(e)). In Lemma \ref{weylem}(e) we saw
$w_{\Delta}w_{\check{J}}\in V^{\check{J}}$. As $V^J\cap V^{\check{J}}=\emptyset$ for $J\ne\check{J}$ this shows $J=\check{J}$ and $w=z^{\check{J}}$, contradicting our hypothesis $w\ne z^J$.\hfill$\Box$\\

The next result concerns the partial ordering $<_{\emptyset}$ of $W$ (i.e. $<_J$ for $J=\emptyset$), called the weak ordering of $W$ in \cite{bjo}.

Assume that the underlying root-system is irreducible and consider the
following subgroup $W_{\Omega}$ of $W$. We write our set of simple roots as
$\Delta=\{\alpha_1,\ldots,\alpha_{l}\}$ and denote by $\alpha_0\in \Phi$ the
unique highest root. Then we define the elements
$\epsilon_1,\ldots,\epsilon_{l}$ in the ${\mathbb R}$-vector space dual to the
one spanned by $\Phi$ by requiring $(\epsilon_i,\alpha_j)=\delta_{ij}$ for
$1\le i,j\le l$. For $1\le i\le l$ we let $w_{\Delta^{(i)}}\in W$ denote the
longest element of the subgroup of $W$ generated by the set
$\{s_{\alpha_j}\,\,|\,\,j\ne
i\}$. Then$$W_{\Omega}-\{1\}=\{w_{\Delta^{(i)}}w_{\Delta}\,\,|\,\,1\le i\le
l,\,\,(\epsilon_i,\alpha_0)=1\}.$$The conjugation action of $W_{\Omega}$ on
$\{s_{\alpha_0},s_{\alpha_1},\ldots,s_{\alpha_{l}}\}$ identifies $W_{\Omega}$
with the automorphism group of the Dynkin diagram of the affine root system
(see \cite{im} pp. 18-20).

\begin{pro}\label{weyllem2} Suppose that the root-system $\Phi$ contains
  no exceptional factor, i.e. that it is a product of root systems of type ${A}$, ${B}$, ${C}$ or ${D}$. There exists a sequence $w_{\Delta}=w_0,w_1,\ldots,w_r=1$ in $W$ such that for all $i\ge1 $ we have $w_{i-1}<_{\emptyset}w_i$, or $w_{i}=uw_{i-1}$ for some $u\in W_{\Omega}$. 
\end{pro}

{\sc Proof:} (I) We first discuss the case where $\Phi$ is irreducible, hence
of type ${A}_l$, ${B}_l$, ${C}_l$ or ${D}_l$ for some $l\in{\mathbb N}$. We use the respective descriptions of $W_{\Omega}$ given in \cite{im} pp. 18-20. We write $s_i=s_{\alpha_i}$.

{\it Case ${A}_{l}$:} Here $W$ can be identified with the symmetric group in
$\{1,\ldots,l+1\}$. We write an element $w\in W$ as the tuple
$[w(1),\ldots,w(l+1)]$. As simple reflections we take the transpositions
$s_i=[1,\ldots,i-1,i+1,i,i+2,\ldots,l+1]\in W$ for $i=1,\ldots,l$. Then
$W_{\Omega}$ consists of the
elements$$w_{\Delta^{(i)}}w_{\Delta}=[i+1,\ldots,l+1,1,\ldots,i]\quad\quad(0\le
i\le l).$$The length $\ell(w)$ of $w\in W$ is the number of all pairs $(i,j)$
with $i<j$ and $w(i)>w(j)$. For $1\le i\le l$ let us
define$$a_i=[l+2-i,\ldots,l+1,l-i+1,\ldots,1],$$$$b_i=[1,\ldots,i,l+1,\ldots,i+1].$$In
particular, $w_{\Delta}=a_1$ and $b_l=1$. Therefore it is enough to show that
for any $1\le i\le l$ we can pass from $a_i$ to $b_i$ by left-multiplication with an
element of $W_{\Omega}$, and that $b_i<_{\emptyset}a_{i+1}$ if $1\le i\le
l-1$. But we indeed have $b_i=w_{\Delta^{(i)}}w_{\Delta}a_i$, whereas, on the other
hand, $b_i<_{\emptyset}a_{i+1}$ follows from $$a_{i+1}=(s_{l-i}\cdots s_1)(s_{l-i+1}\cdots
s_2)\cdots(s_{l-1}\cdots s_i)b_i,$$$$b_i=[1,\ldots,i,l+1,l,\ldots,i+1],$$$$(s_{l-1}\cdots
s_i)b_i=[1,\ldots,i-1,l,l+1,l-1,\ldots,i],$$$$(s_{l-2}\cdots
s_{i-1})(s_{l-1}\cdots
s_i)b_i=[1,\ldots,i-2,l-1,l,l+1,l-2,\ldots,i-1],$$$$(s_{l-3}\cdots
s_{i-1})(s_{l-2}\cdots
s_{i-1})(s_{l-1}\cdots
s_i)b_i=[1,\ldots,i-3,l-2,l-1,l,l+1,l-3,\ldots,i-2]$$etc. from which we see that the
length increases as required.

{\it Case ${B}_{l}$:} Here $W$ can be identified with the group of signed
permutations of $\{\pm 1,\ldots,\pm l\}$, i.e. with all bijections $w:\{\pm
1,\ldots,\pm l\}\to\{\pm 1,\ldots,\pm l\}$ satisfying $-w(a)=w(-a)$ for all
$1\le a\le l$. We write an element $w\in W$ as the tuple
$[w(1),\ldots,w(l)]$. As simple reflections we take the elements
$s_i=[1,\ldots,l-i-1,l-i+1,l-i,l-i+2,\ldots,l]$ for $1\le i\le l-1$, together
with $s_l=[-1,2,\ldots,l]$. Then the length of $w\in W$ can be computed
as $$\ell(w)=|\{\,(i,j)\,\,;\,\,1\le i<j\le l,\,\,w(i)>w(j)\,\}|-\sum_{1\le
  j\le l\atop w(j)<0}w(j)$$(for all this see \cite{bjo} chapter 8.1). The
group $W_{\Omega}$ consists of two elements, its non-trivial element
is$$w_{\Delta^{(1)}}w_{\Delta}=[1,\ldots,l-1,-l].$$For $1\le i\le l$
let$$a_i=[-i,\ldots,-l,i-1,\ldots,1],$$$$b_i=[-i,\ldots,-(l-1),l,i-1,\ldots,1].$$We
pass from $w_{\Delta}$ to $1$ via the
sequence$$w_{\Delta}=[-1,\ldots,-l]=a_1\stackrel{(*)}{\mapsto}b_1<_{\emptyset}a_2\stackrel{(*)}{\mapsto}b_2<_{\emptyset}a_3\stackrel{(*)}{\mapsto}\ldots$$$$\ldots
<_{\emptyset}a_l\stackrel{(*)}{\mapsto}b_l=[l,\ldots,1]\stackrel{(**)}{\mapsto}[1,\ldots,l]=1.$$Here
the relations $b_i<_{\emptyset}a_{i+1}$ result from the equations
$s_{l-i}\cdots s_1 b_i=a_{i+1}$, increasing the length by $l-i$, as one easily
checks. Each step of type $(*)$ is obtained by left-multiplication with
$w_{\Delta^{(1)}}w_{\Delta}$, i.e. $w_{\Delta^{(1)}}w_{\Delta}a_i=b_i$. It
remains to justify the step $(**)$. Observe
that $$w_{\Delta^{(1)}}w_{\Delta}s_{1}\cdots s_l=[l,1,\ldots,l-1].$$Moreover,
for each $w\in W$ satisfying $w(i)>0$ for all $1\le i\le l$ we have
$w<_{\emptyset}s_{1}\ldots s_lw$. Together it follows that, to prove that the
step $(**)$ is permissible, it suffices to show that $(**)$ decomposes into
left-multiplications with (powers of) $[l,1,\ldots,l-1]$ on the one hand, and
with length-increasing left-multiplications with elements of the set
$s_1,\ldots,s_{l-1}$ on the other hand. (Notice that all these operations
preserve the property $w(i)>0$ for all $1\le i\le l$.) But this was shown in
our analysis of case ${A}_{l}$ (or rather $A_{l-1}$), because the $s_1,\ldots,s_{l-1}$ may be viewed
as Coxeter
generators of the symmetric group ${\rm Aut}(\{1,\ldots,l\})$.

{\it Case ${C}_{l}$:} Here $W$ is the same as in case ${B}_{l}$ and we take
the same simple reflections. Again $W_{\Omega}$ consists of two elements, but
this time its non-trivial element
is$$w_{\Delta^{(l)}}w_{\Delta}=[-l,\ldots,-1].$$We pass from $w_{\Delta}$ to
$1$ via the
sequence$$w_{\Delta}=[-1,\ldots,-l]\stackrel{(*)}{\mapsto}[l,\ldots,1]\stackrel{(**)}{\mapsto}[1,\ldots,l]=1.$$Here
$(*)$ is obtained by left-multiplication with $w_{\Delta^{(l)}}w_{\Delta}$. To
justify the step $(**)$ observe
that$$w_{\Delta^{(l)}}w_{\Delta}s_lw_{\Delta^{(l)}}w_{\Delta}s_1\cdots
s_l=[l,1,\ldots,l-1].$$Moreover, for each $w\in W$ satisfying $w(i)>0$ for all
$1\le i\le l$ we have $w<_{\emptyset}s_{1}\cdots s_lw$ (as already noted
above), and $$w_{\Delta^{(l)}}w_{\Delta}s_{1}\cdots
s_lw\quad<_{\emptyset}\quad s_lw_{\Delta^{(l)}}w_{\Delta}s_{1}\cdots s_lw.$$Thus left-multiplication of $[l,1,\ldots,l-1]$ to such $w\in W$ is a permissible operation for our purposes. Therefore we may conclude as in the case ${B}_{l}$.

{\it Case ${D}_{l}$:} Here $W$ can be identified with the group of signed
permutations of $\{\pm 1,\ldots,\pm l\}$ having an even number of negative
entries, i.e. with all bijections $w:\{\pm 1,\ldots,\pm l\}\to\{\pm
1,\ldots,\pm l\}$ satisfying $-w(a)=w(-a)$ for all $1\le a\le l$, and such
that the number $|\{i\,|\,w(i)<0\}|$ is even. We write an element $w\in W$ as
the tuple $[w(1),\ldots,w(l)]$. As simple reflections we take the elements
$s_i$ for $1\le i\le l-1$ used in cases ${B}_{l}$ and ${C}_{l}$, together
with $$s_l=[-2,-1,3,\ldots,l].$$The length of $w\in W$ can be computed (see
\cite{bjo} chapter 8.2)
as $$\ell(w)=|\{\,(i,j)\,\,;\,\,1\le i<j\le l,\,\,w(i)>w(j)\,\}|+|\{\,(i,j)\,\,;\,\,w(i)+w(j)<0\,\}|.$$$W_{\Omega}$
consists of the four elements
$1,w_{\Delta^{(1)}}w_{\Delta},w_{\Delta^{(l-1)}}w_{\Delta}$ and
$w_{\Delta^{(l)}}w_{\Delta}$. Abstractly, if $l$ is even then $W_{\Omega}$ is isomorphic with
${\mathbb Z}/(2)\times{\mathbb Z}/(2)$, with relations
$(w_{\Delta^{(1)}}w_{\Delta})(w_{\Delta^{(l)}}w_{\Delta})=(w_{\Delta^{(l)}}w_{\Delta})(w_{\Delta^{(1)}}w_{\Delta})=w_{\Delta^{(l-1)}}w_{\Delta}$;
if $l$ is odd then $W_{\Omega}$ is isomorphic with
${\mathbb Z}/(4)$, generated by
$w_{\Delta^{(l)}}w_{\Delta}$, with relations
$(w_{\Delta^{(l)}}w_{\Delta})^2=w_{\Delta^{(1)}}w_{\Delta}$ and
$(w_{\Delta^{(l)}}w_{\Delta})^3=w_{\Delta^{(l-1)}}w_{\Delta}$. (We do not need
this.) We have$$w_{\Delta^{(1)}}w_{\Delta}=[-1,2,\ldots,l-1,-l]$$and,
according to the parity of
$l$,$$w_{\Delta^{(l)}}w_{\Delta}=[-l,\ldots,-1]\quad\quad(l\mbox{
  even})$$$$w_{\Delta^{(l)}}w_{\Delta}=[l,1-l,\ldots,-1]\quad\quad(l\mbox{
  odd})$$ (and $w_{\Delta^{(l-1)}}w_{\Delta}=[l,1-l\ldots,-2,1]$ if $l$ is
even, $w_{\Delta^{(l-1)}}w_{\Delta}=[-l,\ldots,-2,1]$ is $l$ is odd). We pass
from $w_{\Delta}$ to $1$ via the
sequence$$w_{\Delta}=[-1,\ldots,-l]\stackrel{(*)}{\mapsto}[l,\ldots,1]\stackrel{(**)}{\mapsto}[1,\ldots,l]=1\quad\quad
(l\mbox{
  even})$$$$w_{\Delta}=[1,-2,\ldots,-l]\stackrel{(*)}{\mapsto}[l,\ldots,1]\stackrel{(**)}{\mapsto}[1,\ldots,l]=1\quad\quad
(l\mbox{ odd}).$$Here $(*)$ is obtained by left-multiplication with
$w_{\Delta^{(l)}}w_{\Delta}$. To justify the step $(**)$ observe
that$$w_{\Delta^{(1)}}w_{\Delta}s_1\cdots s_{l-2}s_l=[l,1,\ldots,l-1].$$For
each $w\in W$ with $w(i)>0$ for all $1\le i\le l-2$ we have
$w<_{\emptyset}s_{1}\cdots s_{l-2}s_lw$. Thus left-multiplication of
$[l,1,\ldots,l-1]$ to such $w\in W$ is a permissible operation for our
purposes and we may conclude as in the case ${B}_{l}$.

(II) In the general case, where $\Phi$ is not necessarily irreducible,
$\Phi$ is a product of root systems as discussed in (I). It is easy to see
that such a product decomposition comes along with a product decomposition of $W$, of $w_{\Delta}$, of $W_{\Omega}$ and of the
ordering $<_{\emptyset}$ (the latter in the obvious sense: $<_{\emptyset}$ is
characterized componentwise). Therefore we may conclude by applying the result
of (I) to all the factors of $\Phi$.\hfill$\Box$\\

\begin{kor}\label{weyllem3} Suppose that the root-system $\Phi$ contains
  no exceptional factor. For each $w\in
  W^J$ there is a sequence $w_0,w_1,\ldots,w_t$ in $W$ (some
  $t\ge0$) with $(w_0)^J=z^J$ and $(w_t)^J=w$ and such that for all $1\le i\le t$ we have $(w_{i})^J=(uw_{i-1})^J$ for some $u\in W_{\Omega}$, or \begin{gather}\ell((w_{i-1})^J)<\ell((w_{i})^J)\quad\mbox{ and }\quad (w_{i})^J=(sw_{i-1})^J\mbox{ for some }s\in S.\label{steigeproperty}\end{gather}
\end{kor}

{\sc Proof:} Observe first that for $w, w'$ in $W$ and
$s\in S$ with $\ell(w')<\ell(w)$ and $w=sw'$ we
have $$[\ell((w')^J)<\ell((w)^J)\,\,\,\,\mbox{
  and }\,\,\,\,(w)^J=s(w')^J=(sw')^J]\quad\quad\mbox{  or  }\quad\quad(w)^J=(w')^J.$$Let $w_{\Delta}=w_0,w_1,\ldots,w_{r}=1$ be a sequence in $W$ such that for all
$1\le i\le r$ we have $w_{i-1}<_{\emptyset}w_i$, or $w_{i}=uw_{i-1}$ for some $u\in
W_{\Omega}$ (Proposition \ref{weyllem2}). We have $(w_0)^J=(w_{\Delta})^J=z^J$ by Lemma
\ref{weylem}(e). By suitably refining the intervals from
$w_{i-1}$ to $w_i$ whenever $w_{i-1}<_{\emptyset}w_i$ we may assume that
whenever $w_{i-1}<_{\emptyset}w_i$ then in addition $w_{i-1}=sw_i$ for some
$s\in S$ (depending on $i$). Then, by the above observation, property (\ref{steigeproperty}) holds true for all
$1\le i\le r$ with $w_{i-1}<_{\emptyset}w_i$; for the other $1\le i\le r$ we
have $(w_{i})^J=(uw_{i-1})^J$ for some $u\in W_{\Omega}$. Choose a reduced expression $w=\sigma_m\cdots\sigma_1$ of
$w$ with $\sigma_i\in S$, then put $t=m+r$ and $w_{i+r}=\sigma_i\cdots\sigma_1$ for $1\le i\le m$. By
the above observation, property (\ref{steigeproperty}) holds true for all $r+1\le i\le t$. We have $w=w_t=(w_t)^J$ since $w\in W^J$.\hfill$\Box$\\

{\bf Remark:} For the irreducible reduced root systems of type $E_8$, $F_4$
and $G_2$ we have $W_{\Omega}=\{1\}$ by \cite{im}. Therefore the statement of
Proposition \ref{weyllem2} cannot hold true in these cases. We do not discuss
the remaining exceptional cases, because we do not know if the statement of Proposition
\ref{weyllem2} holds true for these root systems.

\section{Functions on the Iwahori subgroup}
\label{crisecti}

Let $F$ be a non-Archimedean locally compact field, ${\cal O}_F$ its ring of
integers, $p_F\in{\cal O}_F$ a fixed prime element and $k_F$ its residue
field. Let $G$ be a split connected reductive group over $F$. (Here we commit
the usual abuse of notation: what we really mean is that $G$ is the group of
$F$-rational points of such an algebraic $F$-group scheme, similarly for the
subgroups considered below.) Let $T$ be a split maximal torus, $N\subset G$
its normalizer in $G$ and let $W=N/T$, the corresponding Weyl group. For any
$w\in W$ we choose a representative (with the same name) $w\in N$. Let $P=TU$
be a Borel subgroup with unipotent radical $U$. Let $\Phi\subset X^*(T)={\rm
  Hom}_{alg}(T,{\mathbb G}_m)$ be the set of roots, let $\Phi^+\subset\Phi$ be
the set of $P$-positive roots, let $\Phi^-=\Phi-\Phi^+$, let $\Delta\subset\Phi^+$ be the set of simple roots. Since $T$ is split this root system is reduced.

For $\alpha\in\Phi$ let $U_{\alpha}\subset G$ be the associated root subgroup. Then $U=\prod_{\alpha\in\Phi^+}U_{\alpha}$ (direct product, for any ordering of $\Phi^+$). We need the parabolic subgroups $P_J=PW_JP$ of $G$; each parabolic subgroup of $G$ containing $P$ is of this form (for a suitable $J$). For $w\in W$ let $P_{J,w}=wP_Jw^{-1}$ and let $P^-_{J,w}$ be the parabolic subgroup of $G$ opposite to $P_{J,w}$. We then find$$\Phi-\Phi_J(w)=\{\alpha\in\Phi\,\,|\,\,U_{\alpha}\subset P_{J,w}\}$$or equivalently: $\prod_{\alpha\in\Phi_J(w)}U_{\alpha}$ is the unipotent radical of $P^-_{J,w}$. Note that $P_{J,w}=P_{J,w'}$ for any $w'\in wW_J$.

We choose an Iwahori subgroup $I$ in $G$ compatible with $P$, in the sense that we have the decomposition $$G=\bigcup_{w\in W}IwP$$(disjoint union). For any subgroup $H$ in $G$ we write $H^0=H\cap I$. We will make essential use of the following special case of an important result in the theory of Bruhat and Tits, as recalled in Prop. I.2.2. of \cite{stsc}:   

\begin{pro}\label{brutit} The product map gives a
  bijection $$I=G^0=\prod_{\alpha\in\Phi^+}U_{\alpha}^0\times
  T^0\times\prod_{\alpha\in\Phi^-}U_{\alpha}^0$$for any fixed ordering of
  $\Phi^+$ and $\Phi^-$. 
\end{pro}

\begin{lem} Let $D\subset\Phi$ be a $J$-quasi-parabolic subset. Then $\prod_{\alpha\in D}U_{\alpha}^0$ is a subgroup of $G$ and is independent of the ordering of $D$. We denote it by $U_D^0$.
\end{lem} 

{\sc Proof:} Take any ordering of $D$. Then choose an ordering of $\Phi$ which
restricts to this ordering on $D$ and such that the product
map$$\prod_{\alpha\in \Phi}U_{\alpha}\longrightarrow G$$is injective. Write
$D=\bigcap_{w\in {\Theta}}\Phi_J(w)$ (some ${\Theta}\subset W$). Then of
course $$\prod_{\alpha\in D}U_{\alpha}^0=\bigcap_{w\in
  {\Theta}}\,\prod_{\alpha\in \Phi_J(w)}U_{\alpha}^0$$ (all products
w.r.t. the fixed ordering of $\Phi$, and the intersection is taken inside
$G$). For each $w\in
  {\Theta}$ it follows from Proposition \ref{brutit} that $\prod_{\alpha\in
    \Phi_J(w)}U_{\alpha}^0$ is the intersection of $I$ with the unipotent
  radical of $P^-_{J,w}$. (Notice that Proposition \ref{brutit} holds true for {\it any}
  choice of positive/negative system $(\widetilde{\Phi}^+,\widetilde{\Phi}^-)$ in $\Phi$; here we apply it
  for some $(\widetilde{\Phi}^+,\widetilde{\Phi}^-)$ for which $\Phi_J(w)\subset\widetilde{\Phi}^+$.) In particular, $\prod_{\alpha\in
    \Phi_J(w)}U_{\alpha}^0$ is a subgroup of $G$ and is
  independent of the ordering of $\Phi_J(w)$. Thus, the same statements hold
  true for $\prod_{\alpha\in D}U_{\alpha}^0$ as well.\hfill$\Box$\\

For a topological space ${\mathcal T}$ and an $L$-module $M$ let $C^{\infty}({\mathcal T},M)$ denote the $L$-module of locally constant $M$-valued functions on ${\mathcal T}$.

Applying the functor $C^{\infty}(I,.)$ to the exact sequence (\ref{defmjl}) we
obtain an exact
sequence\begin{gather}C^{\infty}(I,\bigoplus_{{\alpha}\in\Delta-J}L[W^{J\cup\{{\alpha}\}}])\stackrel{}{\longrightarrow}C^{\infty}(I,L[W^{J}])\stackrel{}{\longrightarrow}C^{\infty}(I,{\mathfrak
    M}_J(L)){\longrightarrow}0.\label{grosseq}\end{gather}Observe that we have
natural embeddings, which we view as
inclusions, $$C^{\infty}(I/P^0_{J\cup\{{\alpha}\},w},L)\subset C^{\infty}(I,L),$$$$\bigoplus_{{\alpha}\in\Delta-J\atop w\in
  W^{J\cup\{{\alpha}\}}}C^{\infty}(I/P^0_{J\cup\{{\alpha}\},w},L)\subset \bigoplus_{{\alpha}\in\Delta-J\atop w\in W^{J\cup\{{\alpha}\}}}C^{\infty}(I,L)\cong C^{\infty}(I,\bigoplus_{{\alpha}\in\Delta-J}L[W^{J\cup\{{\alpha}\}}]),$$$$\bigoplus_{w\in W^{J}}C^{\infty}(I/P^0_{J,w},L)\subset C^{\infty}(I,L[W^{J}]),$$by summing over the respective direct summands.

\begin{pro}\label{keyiwa} The sequence$$\bigoplus_{{\alpha}\in\Delta-J\atop w\in W^{J\cup\{{\alpha}\}}}C^{\infty}(I/P^0_{J\cup\{{\alpha}\},w},L)\stackrel{\partial_C}{\longrightarrow}\bigoplus_{w\in W^{J}}C^{\infty}(I/P^0_{J,w},L)\stackrel{\nabla_C}{\longrightarrow}C^{\infty}(I,{\mathfrak M}_J(L))$$obtained by restricting (\ref{grosseq}) is exact.
\end{pro}

{\sc Proof:} {\it Step 1.} We first claim that for any two $J$-parabolic subsets $D$ and
$D'$ of $\Phi$ and for any $\alpha\in\Delta-J$ and $w\in
W^{J\cup\{\alpha\}}(D)$ we have\begin{gather}(U_D^0\cap
  U_{D'}^0)P^0_{J\cup\{\alpha\},w}=(U_D^0 P^0_{J\cup\{\alpha\},w})\bigcap
  (U_{D'}^0 P^0_{J\cup\{\alpha\},w})\label{florian}\end{gather}(where
$AB=(AB)=\{ab\,|\,a\in A, b\in B\}$, but {\it not} (in general) the subgroup generated by $A$ and
$B$). The inclusion $\subset$ is obvious. To prove the inclusion $\supset$ it
is enough to prove\begin{gather}(\prod_{\beta\in D\atop\beta\notin D'}U_{\beta}^0)P^0_{J\cup\{\alpha\},w}\cap
U_{D'}^0\quad\subset\quad P^0_{J\cup\{\alpha\},w}.\label{kosm}\end{gather}Let us
write for the
moment $$\Phi'=\Phi-\Phi_{J\cup\{\alpha\}}(w)=\{\beta\in\Phi\,;\,U_{\beta}\subset
P_{J\cup\{\alpha\},w}\}.$$As $w\in W^{J\cup\{\alpha\}}(D)$ we have $D\cap
\Phi'=\emptyset$. It follows from Proposition \ref{brutit} (applied with a
positive/negative system $(\widetilde{\Phi}^+,\widetilde{\Phi}^-)$ for which
$\widetilde{\Phi}^+\cap \Phi_J(w)$ is before $\widetilde{\Phi}^-\cap \Phi_J(w)$) that we find
subsets ${\mathcal S}_1$ and ${\mathcal S}_2$ of $G^0$ containing the neutral
element, such that$$P^0_{J\cup\{\alpha\},w}=(\prod_{\beta\in
  D'\cap\Phi'}U_{\beta}^0){\mathcal S}_1,$$$$G^0=(\prod_{\beta\in
  D\atop\beta\notin D'}U_{\beta}^0)(\prod_{\beta\in
  D'\cap\Phi'}U_{\beta}^0){\mathcal S}_1(\prod_{\beta\in
  D'\atop\beta\notin\Phi'}U_{\beta}^0){\mathcal S}_2$$and such that all products
are direct (unique factorization of elements). Formula (\ref{kosm}) follows.

{\it Step 2.} Let $(f_w)_{w\in W^J}\in \ke(\nabla_C)$. Choose an enumeration $D_0$, $D_1$, $D_2,\ldots$ of all $J$-quasi-parabolic subsets of $\Phi$ such that $n<m$ implies $|D_n|\le|D_m|$. By induction on $m$ we show: adding to $f$ an element in the image of $\partial_C$ if necessary, we may assume $f_w|_{U_{D_n}^0}=0$ for all $w\in W^J$, all $n\le m$.

Assume we have $f_w|_{U_{D_n}^0}=0$ for all $w\in W^J$, all $n<m$. Let us write $D=D_m$. 

{\it Claim: We have $f_w|_{U_{D}^0}=0$ for all $w\in W^J-W^J(D)$.}

 Indeed, for such $w$ we have $|D\cap\Phi_J(w)|<|D|$, hence
 $D\cap\Phi_J(w)=D_n$ for some
 $n<m$. Thus$$f_w(U_{D}^0)=f_w(U_{D_n}^0\prod_{\alpha\in
   D-D_n}U_{\alpha}^0)=f_w(U_{D_n}^0)=0$$where in the first equation we used
 that we may form $U_{D}^0$ with respect to any ordering of $D$, where the
 second equation follows from $U_{\alpha}^0\subset P^0_{J,w}$ for
 $\alpha\notin \Phi_J(w)$ (and the invariance property of $f_w$), and where
 the last equation holds true by induction hypothesis.

The claim is proven. 

Our sequence in question restricts to a
sequence\begin{gather}\bigoplus_{{\alpha}\in\Delta-J\atop w\in
    W^{J\cup\{{\alpha}\}}(D)}C^{\infty}(I/P^0_{J\cup\{{\alpha}\},w},L)\stackrel{\partial_C^D}{\longrightarrow}\bigoplus_{w\in
    W^{J}(D)}C^{\infty}(I/P^0_{J,w},L)\stackrel{\nabla_C^D}{\longrightarrow}C^{\infty}(I,{\mathfrak
    M}_J(L)).\label{resseq}\end{gather}For any $x\in U_D^0$, evaluating
functions at $x$ transforms (\ref{resseq}) into a sequence isomorphic with the
one from Proposition \ref{parml} (b). Let us denote by $(\partial^D_C)_x$
resp. by $({\nabla^D_C})_x$ the differentials of this sequence, which by
Proposition \ref{parml} (b) is exact. From the above claim it follows that $${f}^D(x)=(f_w(x))_{w\in W^{J}(D)}\in\ke(({\nabla^D_C})_x),$$hence this lies in the image of $(\partial^D_C)_x$. For all $x\in U_D^0$ choose preimages of ${f}^D(x)$ under $(\partial^D_C)_x$. Since the $f_w$ are locally constant, these preimages can be arranged to vary locally constantly on $U_D^0$, and moreover, in view of our induction hypothesis we may assume that for all $x\in U_D^0\cap \cup_{n<m}U^0_{D_n}$ these preimages are zero.

For any ${\alpha}\in\Delta-J$ and $w\in W^{J\cup\{{\alpha}\}}(D)$ the natural
map $U_D^0\to I/P^0_{J\cup\{{\alpha}\},w}$ is injective. Thus we find an
element$$g^D=(g_{\alpha,w})_{\alpha,w}\in\bigoplus_{{\alpha}\in\Delta-J\atop
  w\in
  W^{J\cup\{{\alpha}\}}(D)}C^{\infty}(I/P^0_{J\cup\{{\alpha}\},w},L)$$which on
$U_D^0$ assumes the preimages of the ${f}^D(x)$ just chosen, and which
vanishes at all $x\in\cup_{n<m}U^0_{D_n}$ with $x\notin U_D^0$ --- for this
last property we take advantage of (\ref{florian}). We
obtain$$f^D(x)-\partial_C^D(g^D)(x)=0$$for all $x\in \cup_{n\le m}U^0_{D_n}$:
for $x\in U_{D_m}^0=U_D^0$ this follows from our definition of $g^D|_{U_D^0}$,
for $x\in \cup_{n< m}U^0_{D_n}$ with $x\notin U_D^0$ this follows from the
vanishing of $g^D$ at such $x$ together with the induction hypothesis. Now set
$g_{\alpha,w}=0$ for all ${\alpha}\in\Delta-J$ and $w\in
W^{J\cup\{{\alpha}\}}-W^{J\cup\{{\alpha}\}}(D)$. By the above claim and by what we just saw
we find$$((f_w)_w-\partial_C((g_{\alpha,w})_{\alpha,w}))(x)=0$$for all $x\in
\cup_{n\le m}U^0_{D_n}$. The induction is complete. 

{\it Step 3.} We have shown that, adding to $(f_w)_w\in \ke(\nabla_C)$ an element in the image of $\partial_C$ if necessary, we may assume $f_w|_{U_{D}^0}=0$ for all $w\in W^J$, all $J$-quasi-parabolic subsets $D$. In particular we find $f_w|_{U_{\Phi_J(w)}^0}=0$ for all $w\in W^J$. But $U_{\Phi_J(w)}^0$ is a set of representatives for $I/P^0_{J,w}$ (again invoke Proposition \ref{brutit}), hence $f_w=0$. We are done.\hfill$\Box$\\

{\bf Definition:} Let $J$ be a subset of $\Delta$. We define the
$G$-representation ${\rm Sp}_J(G,L)$ by the exact sequence of
$G$-representations $$\bigoplus_{{\alpha}\in\Delta-J}C^{\infty}(G/P_{J\cup\{\alpha\}},L)\stackrel{\partial}{\longrightarrow}
C^{\infty}(G/P_{J},L)\longrightarrow {\rm Sp}_J(G,L)\longrightarrow0,$$where
$\partial$ is the sum of the canonical inclusions, and the $G$-action is by
left translation of functions on $G$. We call ${\rm Sp}_J(G,L)$ the $J$-special $G$-representation with coefficients in $L$.

\begin{satz}\label{embedding} ${\rm Sp}_J(G,L)$ is $L$-free. There exists an $I$-equivariant embedding$${\rm Sp}_J(G,L)\stackrel{\lambda_L}{\hookrightarrow} C^{\infty}(I,{\mathfrak M}_J(L)).$$Its formation commutes with base changes: for a ring morphism $L\to L'$ the composite$${\rm Sp}_J(G,L)\otimes_LL'\cong{\rm Sp}_J(G,L')\stackrel{\lambda_{L'}}{\hookrightarrow} C^{\infty}(I,{\mathfrak M}_J(L'))\cong C^{\infty}(I,{\mathfrak M}_J(L))\otimes_LL'$$is $\lambda_L\otimes_LL'$.
\end{satz}

{\sc Proof:} Recall that for $w\in W$ we defined $P^0_{J,w}=I\cap
wP_Jw^{-1}$. Note that $P^0_{J,w}$ and $wP_J$ depend only on the coset $wW_J$,
not on the specific representative $w\in wW_J$. The same is true for the
isomorphism $$I/P^0_{J,w}\cong IwP_J/P_J,$$$$i\mapsto iw.$$It follows that for
any inclusion of cosets $wW_J\subset wW_{J\cup\{\alpha\}}$ we have a
commutative diagram$$\xymatrix{I/P^0_{J,w}\ar[d]^{\cong}\ar[r]^{}&
  I/P^0_{J\cup\{\alpha\},w}\ar[d]^{\cong}\\{IwP_J/P_J}\ar[r]^{}&IwP_{J\cup\{\alpha\}}/P_{J\cup\{\alpha\}}}$$where
the horizontal arrows are the obvious projections and the vertical arrows are
the above isomorphisms. Now recall the decompositions$$G/P_J=\cup_{w\in
  W^J}IwP_J/P_J,\quad\quad\quad G/P_{J\cup\{\alpha\}}=\cup_{w\in
  W^{J\cup\{\alpha\}}}IwP_{J\cup\{\alpha\}}/P_{J\cup\{\alpha\}}$$(disjoint
unions). They give$$C^{\infty}(G/P_J,L)=\bigoplus_{w\in
  W^J}C^{\infty}(IwP_J/P_J,L),$$$$C^{\infty}(G/P_{J\cup\{\alpha\}},L)=\bigoplus_{w\in
  W^{J\cup\{\alpha\}}}C^{\infty}(IwP_{J\cup\{\alpha\}}/P_{J\cup\{\alpha\}},L).$$With
these identifications, the above commutative diagrams (for all $\alpha\in
\Delta-J$) induce a commutative diagram$$\xymatrix{\bigoplus_{\alpha\in
    \Delta-J}C^{\infty}(G/P_{J\cup\{{\alpha}\}},L)\ar[d]^{\cong}\ar[r]&
  C^{\infty}(G/P_{J},L)\ar[d]^{\cong}\ar[r]&{\rm
    Sp}_J(G,L)\ar[r]&0\\{\bigoplus}_{{\alpha}\in\Delta-J\atop w\in
    W^{J\cup\{{\alpha}\}}}C^{\infty}(I/P^0_{J\cup\{{\alpha}\},w},L)\ar[r]&\bigoplus_{w\in
    W^{J}}C^{\infty}(I/P^0_{J,w},L)\ar[r]&C^{\infty}(I,{\mathfrak
    M}_J(L))}$$where the vertical arrows are isomorphisms. The top row is
exact by the definition of ${\rm Sp}_J(G,L)$, the bottom row is exact by
Proposition \ref{keyiwa}, and clearly all arrows are $I$-equivariant. Hence we
get the wanted injection $\lambda_L:{\rm Sp}_J(G,L)\hookrightarrow
C^{\infty}(I,{\mathfrak M}_J(L))$. From its construction it is clear that it
commutes with base changes $L\to L'$ as stated. We then derive the freeness of ${\rm Sp}_J(G,L)$: first for $L={\mathbb Z}$ since $C^{\infty}(I,{\mathfrak M}_J({\mathbb Z}))$ is ${\mathbb Z}$-free, then by base change ${\mathbb Z}\to L$ for any $L$.\hfill$\Box$\\

The following corollary was conjectured by Vign\'{e}ras \cite{vig}:

\begin{kor}\label{invber} The submodule ${\rm Sp}_J(G,L)^I$ of $I$-invariants in ${\rm Sp}_J(G,L)$ is free of rank$${\rm rk}_L({\rm Sp}_J(G,L)^I)={\rm rk}_L({\mathfrak M}_J(L))=|V^J|.$$
\end{kor}

{\sc Proof:} By Proposition \ref{parml} we know that ${\mathfrak M}_J(L)$ is
free of rank $|V^J|$. From the definition of ${\mathfrak M}_J(L)$ it follows
that the map$$L[W^J]\cong
\bigoplus_{w\in W^{J}}C^{\infty}(I/P^0_{J,w},L)^I\longrightarrow
C^{\infty}(I,{\mathfrak M}_J(L))^I\cong{\mathfrak M}_J(L)$$is
surjective. In the proof of Theorem \ref{embedding} we saw that the induced
map$${\rm Sp}_J(G,L)^I\longrightarrow
C^{\infty}(I,{\mathfrak M}_J(L))^I\cong{\mathfrak M}_J(L)$$is injective, hence
bijective.\hfill$\Box$\\


\begin{kor} Let $\pi$ be a smooth irreducible (hence finite dimensional) representation of $I$ on a ${\mathbb C}$-vector space. Then $\pi$ occurs in ${\rm Sp}_J(G,{\mathbb C})$ with multiplicity at most $|V^J|\dim_{\mathbb C}(\pi)$. 
\end{kor}

{\sc Proof:} It holds that $\pi$ occurs in $C^{\infty}(I,{\mathfrak M}_J({\mathbb C}))$ with multiplicity $|V^J|\dim_{\mathbb C}(\pi)$.\hfill$\Box$\\

{\bf Remark:} If $L$ is a complete field extension of $F$ we may replace all spaces of locally constant functions occuring here by the corresponding spaces of locally $F$-analytic functions. In particular we may define locally analytic $G$-representations ${\rm Sp}^{\rm an}_J(G,L)$ and $C^{\rm an}(I,{\mathfrak M}_J(L))$. Then Theorem \ref{embedding} and Corollary \ref{invber} carry over, with the same proofs: there exists an $I$-equivariant embedding$${\rm Sp}^{\rm an}_J(G,L)\hookrightarrow C^{\rm an}(I,{\mathfrak M}_J(L))$$and we have ${\rm rk}_L({\rm Sp}^{\rm an}_J(G,L)^I)={\rm rk}_L({\mathfrak M}_J(L))=|V^J|$.

\section{Special representations of finite reductive groups}
 
\label{finredsec}
There is a unique chamber $C$ in the standard apartment associated to $T$ in
the Bruhat-Tits-building of $G$ which is fixed by our Iwahori subgroup
$I$. Let $x_0$ be a special vertex of (the closure of) $C$ and suppose that our Borel subgroup
$P$ is adapted to $x_0$ (see below for what this means). Let ${\mathcal
  G}_{x_0}/{\mathcal O}_F$ denote the ${\mathcal O}_F$-group scheme with
generic fibre the underlying $F$-group scheme ${\mathbb G}$ of $G={\mathbb
  G}(F)$ and such that for each unramified Galois extension $F'$ of $F$ with
ring of integers ${\mathcal O}_{F'}$ we have$${\mathcal G}_{x_0}({\mathcal
  O}_{F'})=\{g\in {\mathbb G}(F')\quad|\quad gx_0=x_0\}$$(see \cite{tit}
section 3.4). This ${\mathcal G}_{x_0}$ is a group scheme as constructed by Chevalley
(\cite{tit} statement 3.4.1). Its special fibre ${\mathcal G}_{x_0}\otimes_{{\mathcal
    O}_{F}}k_F$ is a split connected reductive group over $k_F$ with the same
root datum as $G$ (\cite{tit} statement 3.8.1; compare also \cite{jan}, part II, section 1.17,
and for adjoint semisimple $G$ see \cite{im} p.30/31 where the Bruhat
decomposition of $\overline{G}=({\mathcal G}_{x_0}\otimes_{{\mathcal
    O}_{F}}k_F)(k_F)$ is discussed similarly to how we are going to use it
here). Let $K_{x_0}={\mathcal G}_{x_0}({\mathcal O}_{F})$
and $$U_{x_0}=\ke\quad[\quad K_{x_0}\longrightarrow {\mathcal
  G}_{x_0}(k_{F})\quad].$$For $H$ any of the groups $G$, $P_J$, $P$, $T$, $N$,
$U$, $U_{\alpha}$ let $$\overline{H}=\frac{H\cap K_{x_0}}{H\cap U_{x_0}}.$$Our
requirement above that $P$ be adapted to $x_0$ means that $I$ is the preimage of $\overline{P}$ under the homomorphism $K_{x_0}\to\overline{G}$. On groups of $k_F$-rational points we have: $\overline{P}_J$ is a parabolic subgroup in $\overline{G}$, containing the Borel subgroup $\overline{P}$. This $\overline{P}$ has $\overline{U}$ as its unipotent radical and contains the maximal split torus $\overline{T}$, whose normalizer in $\overline{G}$ is $\overline{N}$. The quotient $\overline{N}/\overline{T}$ is canonically identified with the Weyl group $W=N/T$, and similarly as before we choose for any $w\in W$ a representative (with the same name) $w\in\overline{N}$. Let $\overline{P}^-=\overline{T}\overline{U}^-$ denote the Borel subgroup opposite to $P$, with unipotent radical $\overline{U}^-$. For $w\in W$ let $\overline{U}^w=\overline{U}\cap w\overline{U}^-w^{-1}$. Then$$\overline{U}^w=\prod_{\alpha\in\Phi^+\atop w^{-1}(\alpha)\in\Phi^-}\overline{U}_{\alpha}$$and $\overline{U}^1=\{1\}$. By transposition of \cite{vig} par. 4.2, Prop. 4 (b) we have\begin{gather}\overline{U}^{w}w\overline{P}_J=\overline{P}w\overline{P}_J\label{dirbru}\end{gather}for any $w\in W^J$, and the left hand side product is direct.

\begin{lem}\label{brudec} Let $w\in W^J$ and $s\in S$.\\(a) If $(sw)^J=w$ then $$us\overline{U}^{w}w\overline{P}_J=\overline{U}^{w}w\overline{P}_J$$for each $u\in \overline{U}^s$, and these are direct products.\\(b) If $\ell((sw)^J)>\ell(w)$ then $$\overline{U}^ss\overline{U}^ww\overline{P}_J=\overline{U}^{sw}sw\overline{P}_J$$and these are direct products.\\(c) If $\ell((sw)^J)<\ell(w)$, then $w^{-1}(\beta)\in\Phi^-$, where $s=s_{\beta}$. The product$$\overline{U}'=\prod_{\alpha\in\Phi^+-\{\beta\}\atop w^{-1}(\alpha)\in\Phi^-}\overline{U}_{\alpha}$$(any ordering of the factors) is a subgroup of $\overline{U}^w$. We have$$\overline{U}^ssu\overline{U}'w\overline{P}_J=\overline{U}^{w}w\overline{P}_J\quad\quad\quad\mbox{ for } u\in \overline{U}^{s}-\{1\},$$$$us\overline{U}'w\overline{P}_J=\overline{U}^{sw}sw\overline{P}_J\quad\quad\quad\mbox{ for }u\in \overline{U}^{s}$$and all these are direct products.
\end{lem}

{\sc Proof:} We point out that in all the stated equalities the respective right hand sides are direct products. Therefore, once the equalities are known, the products on the respective left hand sides are seen to be direct simply by a cardinality argument since we work over a finite field.

We use general facts on Bruhat decompositions.\\(a) We
have $$s\overline{U}^ww\overline{P}_J=s\overline{P}w\overline{P}_J\quad\subset\quad\overline{P}w\overline{P}_J
\cup
\overline{P}sw\overline{P}_J=\overline{P}w\overline{P}_J=\overline{U}^ww\overline{P}_J$$where
at the inclusion sign we use $s\overline{P}w\subset \overline{P}w\overline{P}
\cup \overline{P}sw\overline{P}$, and where in the equality following it we
use the hypothesis $(sw)^J=w$, i.e. $swW_J=wW_J$. Applying $s$ we see that
this inclusion is an equality. Since $u\in \overline{P}$ and
$\overline{U}^ww\overline{P}_J=\overline{P}w\overline{P}_J$ we get (a).\\(b)
$\ell((sw)^J)>\ell(w)$ implies $\ell(sw)>\ell(w)$ and again by general
properties of Bruhat decompositions we
find$$\overline{U}^ss\overline{U}^ww\overline{P}_J=\overline{U}^ss\overline{P}w\overline{P}_J=\overline{P}s\overline{P}w\overline{P}_J=\bigcup_{v\in
  W_J}\overline{P}s\overline{P}w\overline{P}v\overline{P}$$$$=\bigcup_{v\in
  W_J}\overline{P}sw\overline{P}v\overline{P}=\overline{P}sw\overline{P}_J=\overline{U}^{sw}sw\overline{P}_J$$where
the assumption $\ell(sw)>\ell(w)$ implied
$\overline{P}s\overline{P}w\overline{P}=\overline{P}sw\overline{P}$, and where
we made repeated use of (\ref{dirbru}) (in the first and in the last equation
with this $J$, and in the second equation by setting $J=\emptyset$ in
(\ref{dirbru})).\\(c) $\ell((sw)^J)<\ell(w)$ implies $\ell(sw)<\ell(w)$ by
Lemma \ref{weylem}(c), hence $w^{-1}(\beta)\in\Phi^-$. One checks that
$\overline{U}'=s\overline{U}^{sw}s$, hence this is a subgroup. Moreover,
$s\overline{U}'=\overline{U}^{sw}s$ and since $\overline{U}^{s}\subset
\overline{P}$ and
$\overline{U}^{sw}sw\overline{P}_J=\overline{P}sw\overline{P}_J$ the last
equality follows. Finally, again by general facts on Bruhat decompositions we
have$$s\overline{U}^ww\overline{P}_J\quad \subset\quad
\overline{U}^ww\overline{P}_J\cup\overline{U}^{sw}sw\overline{P}_J$$and the
union on the right hand side is disjoint (since $swW_J\ne wW_J$). We just saw
that $s\overline{U}'w\overline{P}_J=\overline{U}^{sw}sw\overline{P}_J$, hence
$s(\overline{U}^{w}-\overline{U}')w\overline{P}_J\subset
\overline{U}^ww\overline{P}_J$. It follows
that $$\overline{U}^{s}su\overline{U}'w\overline{P}_J\quad\subset\quad
\overline{U}^{w}w\overline{P}_J$$ for $u\in\overline{U}^{s}-\{1\}$. To see the
reverse inclusion it is enough to show
$\overline{U}'w\overline{P}_J\subset\overline{U}^{s}su\overline{U}'w\overline{P}_J$,
because$$\overline{U}^{s}\overline{U}'=\overline{U}_{\beta}(\prod_{\alpha\in\Phi^+-\{\beta\}\atop w^{-1}(\alpha)\in\Phi^-}\overline{U}_{\alpha})=\prod_{\alpha\in\Phi^+\atop w^{-1}(\alpha)\in\Phi^-}\overline{U}_{\alpha}=\overline{U}^{w}.$$Since $\overline{U}'=s\overline{U}^{sw}s$ this boils down to showing $\overline{U}^{sw}sw\subset s\overline{U}^{s}sus\overline{U}^{sw}sw\overline{P}_J$, i.e. (by (\ref{dirbru})) to $\overline{U}^{sw}sw\subset s\overline{U}^{s}sus\overline{P}sw\overline{P}_J$. A small computation in ${\rm SL}_2(k_F)$ shows that, because of $u\ne1$, there is some $\tilde{u}\in \overline{U}^{s}$ with $s\tilde{u}sus\in \overline{P}$. This implies the wanted inclusion. \hfill$\Box$\\

{\bf Definition:} Similarly as before, we define the $J$-special $\overline{G}$-representation ${\rm Sp}_J(\overline{G},L)$ with coefficients in $L$ by the exact sequence of $\overline{G}$-representations $$\bigoplus_{{\alpha}\in\Delta-J}C(\overline{G}/\overline{P}_{J\cup\{\alpha\}},L)\stackrel{\partial}{\longrightarrow} C(\overline{G}/\overline{P}_{J},L)\longrightarrow {\rm Sp}_J(\overline{G},L)\longrightarrow0.$$Consider the natural map$$C({\overline{G}}/{\overline{P}}_J,L)\longrightarrow C^{\infty}(G/P_J,L),$$$$f\mapsto [g=ky\mapsto f(\overline{k})]$$where we decompose a general element $g\in G$ as $g=ky$ with $k\in K_{x_0}$ and $y\in P_J$ (using the Iwasawa decomposition $G=K_{x_0}P_J$), and where $\overline{k}$ denotes the class of $k$ in ${\overline{G}}=K_{x_0}/U_{x_0}$. We have similar maps for the various ${P}_{J\cup\{\alpha\}}$, hence an embedding\begin{gather}{\rm Sp}_J(\overline{G},L)\hookrightarrow {\rm Sp}_J({G},L).\label{specem}\end{gather}

For the injectivity note e.g. that both sides may be embedded into $C^{\infty}(I,{\mathfrak M}_J(L))$: for the right hand side we saw this in Theorem \ref{embedding}, for the left hand side this can be seen by repeating the construction for $\overline{G}$ instead of $G$.

For $w\in W^J$ we write $$g_w=\chi_{\overline{P}w\overline{P}_J}=\chi_{\overline{U}^ww\overline{P}_J},$$the characteristic function of $\overline{P}w\overline{P}_J=\overline{U}^ww\overline{P}_J$ on $\overline{G}$. We also write $g_w$ for the class of $g_w$ in ${\rm Sp}_J(\overline{G},L)$. 

\begin{pro}\label{ux0in} (a) The embedding (\ref{specem}) induces an isomorphism$${\rm Sp}_J(\overline{G},L)^{\overline{P}}\cong{\rm Sp}_J({G},L)^{I}.$$(b) The set $\{g_w\,\,|\,\,w\in V^J\}$ is an $L$-basis of ${\rm Sp}_J(\overline{G},L)^{\overline{P}}$.
\end{pro}

{\sc Proof:} This follows from Corollary \ref{invber} together with the
$I$-equivariance of the embedding (\ref{specem}). But of course, one could
also directly compute ${\rm Sp}_J(\overline{G},L)^{\overline{P}}$ (i.e. prove
(b)) proceeding as in the proof of Corollary \ref{invber}. Let us also mention
that for $\overline{G}={\rm GL}_n(k_F)$ (some $n$) a proof of (b) is given in
\cite{ss} par.6, and that tor general $\overline{G}$ that proof carries over (this is then similar to \cite{vig} par.4).\hfill$\Box$\\ 

We define the Hecke algebra$${\mathcal H}(\overline{G},\overline{P};L)={\rm End}_{L[{\overline{G}}]}L[{\overline{G}}/{\overline{P}}].$$For a $\overline{G}$-representation on an $L$-vector space $V$ with subspace $V^{\overline{P}}$ of $\overline{P}$-invariants, Frobenius reciprocity tells us that there is an isomorphism$${\rm Hom}_{L[{\overline{G}}]}(L[{\overline{G}}/{\overline{P}}],V)\cong{\rm Hom}_{L[{\overline{P}}]}(L,V)\cong V^{\overline{P}}$$which sends $\psi\in{\rm Hom}_{L[{\overline{G}}]}(L[{\overline{G}}/{\overline{P}}],V)$ to $\psi({\overline{P}})\in V^{\overline{P}}$. Hence $V^{\overline{P}}$ becomes a right ${\mathcal H}(\overline{G},\overline{P};L)$-module. For $g\in\overline{G}$ we define the Hecke operator $T_g\in{\mathcal H}(\overline{G},\overline{P};L)$ by setting$$(T_gf)(h\overline{P})=\sum_{h'\overline{P}\subset h\overline{P}g\overline{P}}f(h'{\overline{P}})$$for $f\in L[{\overline{G}}/{\overline{P}}]$, where for the moment we identify $L[{\overline{G}}/{\overline{P}}]$ with the $L$-module of functions ${\overline{G}}/{\overline{P}}\to L$. For $n\in \overline{N}$ the Hecke operator $T_n$ only depends on the class of $n$ in $W=\overline{N}/\overline{T}$. It acts on $v\in V^{\overline{P}}$ as\begin{gather}vT_n=\sum_{u\in \overline{P}/(\overline{P}\cap n^{-1}\overline{P}n)}un^{-1}v.\label{gbarhec}\end{gather}

Notice that for $s\in S$ we may identify $\overline{U}^s\cong\overline{P}/(\overline{P}\cap s\overline{P}s)$. Thus formula (\ref{gbarhec}) for the Hecke operator $T_s$ acting on $g_w\in {\rm Sp}_J(\overline{G},L)^{\overline{P}}$ becomes\begin{gather}g_wT_{s}=\sum_{u\in\overline{U}^s}(\mbox{the class of }\chi_{us\overline{U}^ww\overline{P}_J})\label{exhec}\end{gather}in ${\rm Sp}_J(\overline{G},L)^{\overline{P}}$.\\

{\it For the rest of this section we assume that $L$ is a field with $\kara(L)=\kara(k_F)$.} 

\begin{lem}\label{weylhec} Let $w\in W^J$ and $s\in S$.\\(a) If $(sw)^J=w$ then $$g_wT_s=0.$$(b) If $\ell((sw)^J)>\ell(w)$ then $$g_wT_s=g_{sw}.$$(c) If $\ell((sw)^J)<\ell(w)$ then $$g_wT_s=-g_w.$$
\end{lem}

{\sc Proof:} This follows from Lemma \ref{brudec} and from
$|\overline{U}^s|=0$ in $L$. For example, for (c) we compute, using the
notations of Lemma \ref{brudec} (c), in particular the direct product
decomposition
$\overline{U}^w=\overline{U}^s\overline{U}'$: $$g_wT_{s}=\sum_{u\in\overline{U}^s}[\chi_{us\overline{U}^ww\overline{P}_J}]=\sum_{u\in\overline{U}^s}\sum_{u'\in\overline{U}^s}[\chi_{usu'\overline{U}'w\overline{P}_J}]$$$$=\sum_{u\in\overline{U}^s}\sum_{u'\in\overline{U}^s-\{1\}}[\chi_{usu'\overline{U}'w\overline{P}_J}]+\sum_{u\in\overline{U}^s}[\chi_{us\overline{U}'w\overline{P}_J}].$$Lemma
\ref{brudec} (c) together with $|\overline{U}^s|=0$ in $L$ shows that the
second term vanishes and that the first term is
$-[\chi_{\overline{U}w\overline{P}_J}]$. For statement (b) notice that by
Lemma \ref{weylem}(b) we have $sw\in W^J$ (and even $sw\in V^J$ if $w\in
V^J$). (Of course, Lemma \ref{weylhec} may also be deduced from general facts on
Iwahori Hecke algebras; we have included the proof in order to keep the
presentation self contained.)\hfill$\Box$

\begin{pro}\label{indeco} Each non-zero ${\mathcal H}(\overline{G},\overline{P};L)$-submodule $E$ of ${\rm Sp}_J(\overline{G},L)^{\overline{P}}$ contains the element $g_{z^J}$. In particular, the ${\mathcal H}(\overline{G},\overline{P};L)$-module ${\rm Sp}_J(\overline{G},L)^{\overline{P}}$ is indecomposable.   
\end{pro}

{\sc Proof:} Choose an enumeration $z^J=w_0, w_1, w_2,\ldots$ of $V^J$ such
that $w_j<_Jw_i$ implies $i<j$. By Proposition \ref{ux0in} we may write any
element $h$ of $E$ as$$h=\sum_{w\in V^J}\beta_w(h)g_w$$ with certain uniquely determined
$\beta_w(h)\in L$. For $t\ge0$ define the subset $${\mathfrak P}(t)\quad=\quad\{\,\,h\in
   E\,\,|\,\,\beta_{w_i}(h)=0\mbox{ for all }i>t\,\,\mbox{ and
   }\beta_{w_t}(h)\ne0\,\,\}$$of $E$. It is enough to show ${\mathfrak
     P}(0)\ne\emptyset$. As $E-\{0\}=\cup_{t\ge0}{\mathfrak P}(t)$ it is
   enough to show the following: If ${\mathfrak P}(t)\ne\emptyset$ for some
   $t>0$, then ${\mathfrak P}(t')\ne\emptyset$ for some $0\le
 {t}'<t$. 

 By Lemma \ref{weyllem1}, applied to $w_t\in V^{J}-\{z^J\}$, we find some
 $w'\in V^J$ and some $s\in S$
 with$$w_t<_Jw',\quad\quad\ell((sw_t)^J)<\ell(w_t),\quad\quad\ell((sw')^J)\ge\ell(w').$$By
 the definition of $w_t<_Jw'$ we find $s_1,\ldots s_r\in S$ such that, setting
 $w^{(g)}=(s_g\cdots s_1w_t)^J$ for $0\le g\le r$, we have$$\ell(w^{(g)})<\ell(w^{(g+1)})\mbox{ for all }0\le g\le r,\quad\quad\mbox{ and
}\quad\quad w^{(r)}=w'.$$From Lemma \ref{weylem}(f) it follows that in fact
$w^{(i)}\in V^J$ for all $i$. Since we have
$\ell((sw^{(r)})^J)\ge\ell(w^{(r)})$, a case by case inspection of Lemma
\ref{weylhec} shows that $\beta_{w^{(r)}}(E\cdot T_s)=0$. We pick some $h\in
{\mathfrak P}(t)$ and make the following 

{\it Claim: We have $h T_s\in{\mathfrak P}(t)$ and $\beta_{w^{(r)}}(hT_s)=0$.} 

By what we just said, we have $\beta_{w^{(r)}}(hT_s)=0$. Next, we have $hT_s\in\{0\}\cup(\cup_{{t}'\le t}{\mathfrak
    P}({t}'))$ as follows from Lemma \ref{weylhec}, again a case by case
  inspection. Thus it remains to show $\beta_{w_t}(hT_s)\ne 0$. From
  $\ell((sw_t)^J)<\ell(w_t)$ we deduce, again using Lemma \ref{weylhec}, that
  $\beta_{w_t}(hT_s)=-\beta_{w_t}(h)+\beta_{(sw_t)^J}(h)$ if $(sw_t)^J\in
  V^J$, but $\beta_{w_t}(hT_s)=-\beta_{w_t}(h)$ if $(sw_t)^J\notin
  V^J$. On the other hand, if $(sw_t)^J\in
  V^J$ then from $\ell((sw_t)^J)<\ell(w_t)$ we also deduce
  $\beta_{(sw_t)^J}(h)=0$ since $h\in
{\mathfrak P}(t)$. In either case we get $\beta_{w_t}(hT_s)=-\beta_{w_t}(h)\ne
0$. The claim is proven.

In view of this claim we see that there is some $h\in
{\mathfrak P}(t)$ with $\beta_{w^{(r)}}(h)=0$. 

{\it Claim: At least one of the following statements hold true: (a) ${\mathfrak P}({t}')\ne\emptyset$ for some $0\le
 {t}'<t$, or (b) for any $1\le g\le r$ there is some $h\in
{\mathfrak P}({t})$ with $\beta_{w^{(g)}}(h)=0$.}

Assume that (a) is false. Then we prove (b) by descending induction on $g$. For $g=r$ this was just
done. Now let $1\le g<r$ and let $h\in
{\mathfrak P}(t)$ be such that $\beta_{w^{(g+1)}}(h)= 0$. If also
$\beta_{w^{(g)}}(h)=0$ then we are done for this $g$, thus we assume $\beta_{w^{(g)}}(h)\ne 0$. 

 Since we have $\ell(w^{(g)})<\ell(w^{(g+1)})$, Lemma \ref{weylhec} shows
$$\beta_{w^{(g+1)}}(hT_{s_{g+1}})=\beta_{w^{(g)}}(h)\quad\quad\mbox{ and }\quad\quad\beta_{w^{(g)}}(hT_{s_{g+1}})=0.$$As argued similarly in the previous claim,
Lemma \ref{weylhec} also shows $hT_{s_{g+1}}\in\{0\}\cup(\cup_{{t}'\le t}{\mathfrak
    P}({t}'))$. But $hT_{s_{g+1}}\ne0$ since
  $\beta_{w^{(g+1)}}(hT_{s_{g+1}})=\beta_{w^{(g)}}(h)\ne0$, thus $hT_{s_{g+1}}\in{\mathfrak
    P}({t}')$ for some $0\le {t}'\le t$. As we assume that (a) is
  false this means $hT_{s_{g+1}}\in{\mathfrak
    P}({t})$. The claim is proven.

Of course, the last argument applies in the same way for $g=0$: but
since there is no $h\in
{\mathfrak P}({t})$ with $\beta_{w^{(0)}}(h)=\beta_{w_t}(h)=0$, the result is that
indeed ${\mathfrak P}({t}')\ne\emptyset$ for some $0\le{t}'<t$. We are done.\hfill$\Box$\\ 

\begin{kor}\label{pwdiff} The ${\mathcal H}(\overline{G},\overline{P};L)$-modules ${\rm Sp}_J(\overline{G},L)^{\overline{P}}$ for different $J\subset\Delta$ are pairwise non-isomorphic. 
\end{kor}

{\sc Proof:} (That this follows from Proposition \ref{indeco} and Lemma
\ref{weylhec} was pointed out to me by Florian Herzig.) It follows from
Proposition \ref{indeco} that ${\rm Sp}_J(\overline{G},L)^{\overline{P}}$
contains a unique irreducible ${\mathcal
  H}(\overline{G},\overline{P};L)$-submodule ${\mathcal M}_J$. Like any
irreducible ${\mathcal H}(\overline{G},\overline{P};L)$-module it must be
one-dimensional. Therefore Lemma \ref{weylhec} together with Proposition
\ref{indeco} show that $T_s$ for $s\in S$ acts on ${\mathcal M}_J$ with
eigenvalues $0$ or $-1$. More precisly, $T_s$ acts with eigenvalue $0$ if
$(sz^J)^J=z^J$, and with eigenvalue $-1$ if $\ell((sz^J)^J)<\ell(z^J)$, and by
Lemma \ref{weylem} no other cases occur. In fact, Lemma \ref{weylem} says that
$(sz^J)^J=z^J$ is equivalent with $\ell(sz^J)>\ell(z^J)$, and
$\ell((sz^J)^J)<\ell(z^J)$ is equivalent with
$\ell(sz^J)<\ell(z^J)$. Thus$$\{s\in S\,|\,T_s|_{{\mathcal M}_J}=0\}=\{s\in
S\,|\,\ell(sz^J)<\ell(z^J)\},$$but this set allows us to recover $J$. Indeed, let $\check{J}=-w_{\Delta}(J)\subset \Delta$, or equivalently, $\check{J}$ is the subset of $\Delta$ with $w_{\check{J}}=w_{\Delta}w_Jw_{\Delta}$ and $w_{J}=w_{\Delta}w_{\check{J}}w_{\Delta}$. Then $w_{\Delta}=w_{\check{J}}z^J$ (as $z^J=w_{\Delta}w_J$), and since $\ell(w_{\Delta})=\ell(w_{\check{J}})+\ell(z^J)$ we see that $\ell(sz^J)<\ell(z^J)$ for $s\in S$ is equivalent with $\ell(w_{\check{J}}s)>\ell(w_{\check{J}})$, and this is equivalent with $s\notin \check{J}$. But ${J}=-w_{\Delta}(\check{J})$.  \hfill$\Box$\\      

\section{Irreducibility in the residual characteristic}

Now assume for simplicity that $G$ is semisimple. Following our conventions we
put $T^0=I\cap T$ and then let $\widetilde{W}=N/T^0$. This group acts on the apartment $A$ and
can be canonically identified with the semidirect product $(T/T^0)\rtimes
W$. (The embedding $W\to\widetilde{W}$ sends an element of $W=N(T)/T$ to its
unique representative in $\widetilde{W}=N/T^0$ which fixes $x_0$.) It
contains the affine Weyl-group $W^a$, the subgroup of $\widetilde{W}$
generated by the reflections in the walls of $A$. On the other hand, let
$\Omega$ be the subgroup of $\widetilde{W}$ stabilizing the standard chamber
in $A$ (i.e. the one fixed by $I$). Then $\widetilde{W}$ is canonically
identified with the semidirect product $W^a\rtimes\Omega$. If $G$ is of adjoint type the canonical projection $\varphi:\widetilde{W}\to W$ is injective on $\Omega$ and its image $W_{\Omega}=\varphi(\Omega)\subset W$ coincides with the one defined in section \ref{weysec}.\\

We define the Iwahori Hecke algebra$${\mathcal H}(G,I;L)={\rm End}_{L[G]}L[G/I].$$For a smooth $G$-representation on an $L$-vector space $V$ with subspace $V^I$ of $I$-invariants, Frobenius reciprocity tells us that there is an isomorphism$${\rm Hom}_{L[G]}(L[G/I],V)\cong{\rm Hom}_{L[I]}(L,V)\cong V^I$$which sends $\psi\in{\rm Hom}_{L[G]}(L[G/I],V)$ to $\psi(I)\in V^I$. Hence $V^I$ becomes a right ${\mathcal H}(G,I;L)$-module. For $g\in G$ we define the Hecke operator $T_g\in{\mathcal H}(G,I;L)$ by setting$$(T_gf)(hI)=\sum_{h'I\subset hIgI}f(h'I)$$for $f\in L[G/I]$, where for the moment we identify $L[G/I]$ with the $L$-module of compactly supported functions ${G}/I\to L$. The Hecke operator $T_n$ for $n\in N$ depends only on the class of $n$ in $\widetilde{W}$, and the $T_n$ for $n$ running through a system of representatives for $\widetilde{W}$ form an $L$-basis of ${\mathcal H}(G,I;L)$ (\cite{vihec} section 1.3, example 1). They act on $v\in V^I$ as$$vT_n=\sum_{u\in I/(I\cap n^{-1}In)}un^{-1}v.$$

By  Proposition \ref{ux0in} we have an isomorphism\begin{gather}{\rm
    Sp}_J(\overline{G},L)^{\overline{P}}\cong {\rm
    Sp}_J({G},L)^I.\label{invcanucap}\end{gather}For $w\in W$ we had defined a
Hecke operator $T_w$ acting on the ${\mathcal
  H}(\overline{G},\overline{P};L)$-module ${\rm
  Sp}_J(\overline{G},L)^{\overline{P}}$. On the other hand, if we denote again
by $w$ a representative in $N$ of the image of $w$ in $\widetilde{W}$ (under
the embedding $W\hookrightarrow(T/T^0)\rtimes W\cong\widetilde{W}$), we get a
Hecke operator $T_w$ acting on the ${\mathcal H}(G,I;L)$-module ${\rm
  Sp}_J({G},L)^I$. (Note however that, for fixed Iwahori subgroup $I$, the
isomorphism $(T/T^0)\rtimes W\cong\widetilde{W}$ and hence the embedding $W\to
\widetilde{W}$ depends on the choice of the special vertex $x_0$ in (the
closure of) the chamber $C$ fixed by $I$. Hence the ${\mathcal H}(G,I;L)$-elements $T_w$ for $w\in W$ depend on this choice.) It is clear from our constructions that these actions coincide under our isomorphism (\ref{invcanucap}). Recall that for $w\in W^J$ we wrote $g_w$ for the class in ${\rm Sp}_J(\overline{G},L)^{\overline{P}}$ of the characteristic function of $\overline{P}w\overline{P}_J$ on $\overline{G}$. Now we also write $g_w$ for its image in ${\rm Sp}_J({G},L)^I$ under (\ref{invcanucap}), i.e. for the class in ${\rm Sp}_J({G},L)^I$ of the characteristic function of $IwP_J$ on $G$.\\

{\it For the rest of this section we assume that $L$ is a field with $\kara(L)=\kara(k_F)$.} 

\begin{lem}\label{liftome} Assume that $G$ is of adjoint type. For each $u\in W_{\Omega}$ there exists a lifting $\widetilde{u}\in N$ (under the canonical projections $N\to \widetilde{W}\to W$) which normalizes $I$ and such that for all $w\in W^J$ we have $g_wT_{\widetilde{u}^{-1}}=g_{(uw)^J}$ in ${\rm Sp}_J({G},L)^I$. 
\end{lem}

{\sc Proof:} By \cite{im} Proposition 2.10 we can lift $u\in W_{\Omega}$ to an
element $\widetilde{u}\in N$ which normalizes $I$. Therefore
$T_{\widetilde{u}^{-1}}$ acts on ${\rm Sp}_J({G},L)^I$ simply through the
action of $\widetilde{u}\in N\subset G$ and for $w\in W^J$ we compute
$\widetilde{u}IwP_J=I{\widetilde{u}w}P_J=I({uw})^JP_J$. The Lemma
follows.\hfill$\Box$\\

 The hypothesis that $G$ be of adjoint type should be superfluous for Lemma
 \ref{liftome} (if $W_{\Omega}$ is replaced with $\varphi(\Omega)$), but
 \cite{im} assumes this. However, the proof of Theorem \ref{hecirr} below forces us to pass to
the adjoint quotient of $G$ anyway, i.e. for a more serious reason.

\begin{satz}\label{hecirr} If the root-system $\Phi$ contains
  no exceptional factor then the ${\mathcal H}(G,I;L)$-module ${\rm Sp}_J({G},L)^I$ is irreducible.
\end{satz}

{\sc Proof:} By Proposition \ref{indeco} we know that each non-zero ${\mathcal
  H}({G},{I};L)$-submodule of ${\rm Sp}_J({G},L)^{{I}}$ contains the element
$g_{z^J}$. Therefore it is enough to show that ${\rm Sp}_J({G},L)^{{I}}$ is
generated as a ${\mathcal H}({G},{I};L)$-module by the element $g_{z^J}$.\\(a)
We first assume that $G$ is of adjoint type. We claim that for each subspace
$E$ of ${\rm Sp}_J({G},L)^I$ containing $g_{z^J}$ and stable under all $T_w$
for $w\in W$, and stable under all $T_{\widetilde{u}^{-1}}$ for
$\widetilde{u}\in N$ normalizing $I$ as in Lemma \ref{liftome}, we have
$E={\rm Sp}_J({G},L)^I$. Indeed, we know that ${\rm Sp}_J({G},L)^{{I}}$ is
generated as an $L$-vector space by all $g_w$ for $w\in V^J$, so we need to prove $g_{w}\in E$ for each such $w\in
V^J$. To do this we choose a sequence $w_0,w_1,\ldots,w_t$ in $W$  with $(w_0)^J=z^J$ and $(w_t)^J=w$ and such that for all $i\ge1$
we have $(w_{i}) ^J=(uw_{i-1})^J$ for some $u\in W_{\Omega}$,
or $$\ell((w_{i-1})^J)<\ell((w_{i})^J)\quad\mbox{ and
}\quad(w_{i})^J=(sw_{i-1})^J\mbox{ for some }s\in S.$$Such a sequence does exist as we learn from Corollary \ref{weyllem3} because, since we assume that $G$ is of adjoint
type, we may lift the elements of $W_{\Omega}$ to elements of $N$. Now we use Lemmata
\ref{liftome} and \ref{weylhec}(b) to prove by induction on $i$ that $g_{(w_i)^J}\in E$ for all $0\le i\le t$: for $i=0$ this is the hypothesis $g_{z^J}\in E$, for $i=t$ this is the statement $g_{w}\in E$ which we needed to prove.

(b) In the general case we find a central isogeny $\pi:G\to G'$ with $G'$
split, connected, semisimple and of adjoint type, and with the same root
system. We find a split maximal torus $T'$ with normalizer $N'$, a Borel
subgroup $P'$ and an Iwahori subgroup $I'$ in $G'$ such that $\pi^{-1}(T')=T$,
$\pi^{-1}(P')=P$, $\pi^{-1}(I')=I$ and such that $W\cong N'/T'$ (observe that $G$
is semisimple, hence its finite center is contained in $I$). As $\ker(\pi)\subset T$ it is clear that $\pi$ induces a $G$-equivariant isomorphism ${\rm Sp}_J({G'},L)\cong{\rm Sp}_J({G},L)$ which restricts to an isomorphism of Iwahori invariant spaces ${\rm Sp}_J({G'},L)^{I'}\cong{\rm Sp}_J({G},L)^I$ (both of dimension $|V^J|$, by Corollary \ref{invber}). 

We identify the Bruhat-Tits buildings of $G$ and $G'$; then $C$ is fixed by
$I'$, and $P'\subset G'$ is adapted to $x_0$. Let $\widetilde{u}\in N'$ as in Lemma \ref{liftome}, in
particular normalizing $I'$. For $n'\in N'$ we
have\begin{gather}T_{n'}T_{\widetilde{u}^{-1}}=T_{n'{\widetilde{u}^{-1}}}=T_{\widetilde{u}^{-1}}T_{\widetilde{u}n'{\widetilde{u}^{-1}}}\quad\quad\mbox{in
  }{\mathcal H}(G',I';L)\label{braid}\end{gather}by general facts on
${\mathcal H}(G',I';L)$ (the 'braid relations'), or just by the definition of
the $T_g$'s. Now $\widetilde{u}\pi(N)\widetilde{u}^{-1}=\pi(N)$ because $\pi$
is a central isogeny, and this is
contained in $N'$. Since ${\mathcal H}(G,I;L)$ is generated by the $T_n$ with
$n\in N$ (see, e.g. \cite{vihec} section 1.3, example 1), the relations
(\ref{braid}) imply\begin{gather}{\mathcal
    H}(G,I;L)T_{\widetilde{u}^{-1}}=T_{\widetilde{u}^{-1}}{\mathcal
    H}(G,I;L)\label{twinhec}\end{gather}inside ${\rm End}_L({\rm
  Sp}_J({G},L)^I)^{\rm op}$ (here we keep the names of ${\mathcal H}(G,I;L)$ and $T_{\widetilde{u}^{-1}}$ also for their images in ${\rm End}_L({\rm Sp}_J({G},L)^I)^{\rm op}$). We get\begin{gather}(g_{z^J}{\mathcal H}(G,I;L))T_{\widetilde{u}^{-1}}\quad\subset\quad ({\widetilde{u}}g_{z^J}){\mathcal H}(G,I;L)\label{incl1}\end{gather}inside ${\rm Sp}_J({G},L)^I$ (recall that $T_{\widetilde{u}^{-1}}$ acts from the right on ${\rm Sp}_J({G},L)^I$ by left multiplication with $\widetilde{u}$). By Proposition \ref{indeco} we have $g_{z^J}\in({\widetilde{u}^{-1}}g_{z^J}){\mathcal H}(G,I;L)$. We apply $T_{\widetilde{u}^{-1}}$, by equation (\ref{twinhec}) again this gives ${\widetilde{u}}g_{z^J}\in g_{z^J}{\mathcal H}(G,I;L)$, and together with (\ref{incl1}) we get$$(g_{z^J}{\mathcal H}(G,I;L))T_{\widetilde{u}^{-1}}\quad\subset\quad g_{z^J}{\mathcal H}(G,I;L).$$By what we have seen in (a) this proves the Theorem.\hfill$\Box$\\

{\bf Remarks:} (a) We just saw that, in case $\Phi$ contains
  no exceptional factor (possibly also factors $E_6$, $E_7$ can be allowed,
  see the remark at the end of section \ref{weysec}), to prove the
  irreducibility of the ${\mathcal H}(G,I;L)$-module ${\rm Sp}_J({G},L)^{{I}}$
  it is enough to use the action of ${\mathcal
    H}(\overline{G},\overline{P};L)$ together with the Hecke operators
  $T_{\widetilde{u}^{-1}}$ of Lemma \ref{liftome}. 

(b) Corollary \ref{invber} together with \cite{vig} Proposition 10 provides us with an isomorphism of ${\mathcal H}(G,I;L)$-modules \begin{gather}{\rm Sp}_J({G},L)^{{I}}\cong\frac{C^{\infty}(G/P_J,L)^I}{\sum_{\alpha\in\Delta-J}C^{\infty}(G/P_{J\cup\{\alpha\}},L)^I}.\label{rachel}\end{gather}

\begin{kor}\label{girr} If the root-system $\Phi$ contains
  no exceptional factor then the $G$-representation ${\rm Sp}_J({G},L)$ is irreducible.
\end{kor}

{\sc Proof:} Let $I_1\subset I$ denote the pro-$p$-Iwahori subgroup in $I$, where $p=\kara(k_F)$. Then $I$ is generated by $I_1$ and $T^0=T\cap I$. By Proposition \ref{keyiwa} and the proof of Theorem \ref{embedding} we may identify ${\rm Sp}_J({G},L)$ as an $L[I]$-module with the image of ${\nabla_C}$ (notation of Proposition \ref{keyiwa}). As such it is contained in $C^{\infty}(I/T^0,{\mathfrak M}_J(L))$. Since we obviously have $C^{\infty}(I/T^0,{\mathfrak M}_J(L))^{I_1}=C^{\infty}(I/T^0,{\mathfrak M}_J(L))^{I}$ it follows that $${\rm Sp}_J({G},L)^I={\rm Sp}_J({G},L)^{I_1}.$$(This argument was suggested by Vign\'{e}ras.) Replacing $I$ by $I_1$ in our definition of the Iwahori Hecke Algebra ${\mathcal H}(G,I;L)$ we obtain the algebra ${\mathcal H}(G,I_1;L)$. Similarly as before, ${\rm Sp}_J({G},L)^{I_1}$ is an ${\mathcal H}(G,I_1;L)$-module, and the irreducibility of ${\rm Sp}_J({G},L)^I$ as an ${\mathcal H}(G,I;L)$-module (Theorem \ref{hecirr}) immediately implies the irreducibility of ${\rm Sp}_J({G},L)^{I_1}={\rm Sp}_J({G},L)^{I}$ as an ${\mathcal H}(G,I_1;L)$ module. Now recall the well known fact that for every smooth representation of a pro-$p$-group --- like $I_1$ --- on a non-zero $L$-vector space $E$ the subspace $E^{I_1}$ of $I_1$-invariants is non-zero (since $\kara(L)=p$). Applied to a non-zero $G$-subrepresentation $E$ of ${\rm Sp}_J({G},L)$, the irreducibility of ${\rm Sp}_J({G},L)^{I_1}$ as a ${\mathcal H}(G,I_1;L)$ module implies $E^{I_1}={\rm Sp}_J({G},L)^{I_1}$. But ${\rm Sp}_J({G},L)$ is generated as a $L[G]$-module by ${\rm Sp}_J({G},L)^{I_1}$; this follows from \cite{vig}, Proposition 9, where it is shown that even the $L[G]$-module $C^{\infty}(G/P_J,L)$ is generated by its $I_1$-fixed vectors. Thus $E={\rm Sp}_J({G},L)$ and we are done.\hfill$\Box$\\

{\bf Remarks:} (a) For any $J$ with $|V^{J}|=1$, like $J=\emptyset$, we get the
irreducibility of ${\rm Sp}_J({G},L)$ for {\it any} $G$ (even if $\Phi$
contains exceptional factors). The irreducibility of the
Steinberg representation ${\rm Sp}_{\emptyset}({G},L)$ had been obtained
earlier by Vign\'{e}ras \cite{vig}. See \cite{her} for the irreducibility
statement in general.

(b) Vign\'{e}ras \cite{vig} shows that each ${\rm
  Sp}_J({G},L)$ admits a $P$-equivariant filtration, with factors the natural
$P$-representations $C^{\infty}_c(PwP/P,L)$ for $w\in V^J$. These factors are
shown to be irreducible (\cite{vig} Proposition 1, Theorem 5).

\begin{kor}\label{johoe} (a) The $G$-representations ${\rm Sp}_J({G},L)$ for
  the various subsets $J\subset \Delta$ are pairwise non-isomorphic.\\(b)
  Suppose that the root-system $\Phi$ contains
  no exceptional factor. The $G$-representations ${\rm Sp}_J({G},L)$ with $J$ running through all subsets $J\subset \Delta$ form the irreducible constituents of the $G$-representation $C^{\infty}(G/P,L)$, each one occuring with multiplicity one. 
\end{kor}

{\sc Proof:} Statement (a) follows from Corollary \ref{pwdiff}. The
irreducibility of the ${\rm Sp}_J({G},L)$ in (b) is Corollary \ref{girr}. Now put
$F_{-1}=0\subset C^{\infty}(G/P,L)$ and$$F_i=\sum_{J\subset\Delta\atop|J|=|\Delta|-i}C^{\infty}(G/P_J,L)$$for
$i\ge0$. Then $0=F_{-1}\subset F_0\subset F_1\subset\ldots\subset F_{|\Delta|}=C^{\infty}(G/P,L)$ is an exhaustive $G$-equivariant
filtration. To prove the remaining statements in (b) it is enough to see
that for any $i\ge0$ there exists a $G$-equivariant
isomorphism\begin{gather}\frac{F_i}{F_{i-1}}\cong\bigoplus_{J\subset\Delta\atop|J|=|\Delta|-i}{\rm
  Sp}_J(G,L).\label{floriansfrage}\end{gather}We do this by induction on $i$. For any $J\subset\Delta$ with
$|J|=|\Delta|-i$ we have a natural $G$-equivariant
map $C^{\infty}(G/P_J,L)\to F_i$, inducing an embedding $$\iota_J:{\rm
  Sp}_J({G},L)\hookrightarrow\frac{F_i}{\sum_{\alpha\in \Delta-J}C^{\infty}(G/P_{J\cup\{\alpha\}},L)}.$$From the induction hypothesis, from the
  irreducibility of the ${\rm Sp}_{J'}({G},L)$ and their being pairwise
  non-isomorphic it follows that $\iota_J$ induces an embedding ${\rm
  Sp}_J({G},L)\hookrightarrow{F_i}/F_{i-1}$. Next, from the
  irreducibility of the ${\rm Sp}_J({G},L)$ and their being pairwise
  non-isomorphic again, it follows that these embeddings sum up to an
  isomorphism (\ref{floriansfrage}) as desired.\hfill$\Box$\\ 

{\bf Question:} Is the theory of extensions between the various $G$-representations ${\rm Sp}_J({G},L)$ (for $L$ a field with $\kara(L)=\kara(k_F)$) parallel to the theory of extensions between the various $G$-representations ${\rm Sp}_J({G},{\mathbb{C}})$ (as worked out in \cite{orl}, \cite{ss}) ?

\begin{kor}\label{unilat} Suppose that the root-system $\Phi$ contains
  no exceptional factor. Let ${\mathcal O}_K$ be a complete discrete valuation ring with fraction field $K$ and residue field $k_K$. Suppose $\kara(k_K)=\kara(k_F)$. Up to $K^{\times}$-homothety, ${\rm Sp}_J({G},{\mathcal O}_K)$ is the unique $G$-stable ${\mathcal O}_K$-lattice inside ${\rm Sp}_J({G},K)$.
\end{kor}

{\sc Proof:} (I thank Marie-France Vign\'{e}ras for completing my argument
here.) Let ${\mathcal S}$ be another $G$-stable ${\mathcal O}_K$-lattice inside ${\rm
  Sp}_J({G},K)$. Let $p_K\in {\mathcal O}_K$ be a uniformizer. Since ${\rm
  Sp}_J({G},k_K)$ is irreducible by Corollary \ref{girr}, the image of
$p_K^n{\mathcal S}\cap {\rm Sp}_J({G},{\mathcal O}_K)$ in ${\rm Sp}_J({G},{\mathcal
  O}_K)\otimes_{{\mathcal O}_K} k_K={\rm Sp}_J({G},k_K)$ for $n\in{\mathbb Z}$
must be either (a) zero, or (b) all of ${\rm Sp}_J({G},k_K)$. Case (a) implies
$p_K^{n-1}{\mathcal S}\subset{\rm Sp}_J({G},{\mathcal O}_K)$. Case (b)
implies \begin{gather}{\rm Sp}_J({G},{\mathcal O}_K)\quad\subset\quad p_K{\rm
    Sp}_J({G},{\mathcal O}_K)+p_K^n{\mathcal S}.\label{infini}\end{gather}Now ${\rm
  Sp}_J({G},{\mathcal O}_K)$ is finitely generated as an ${\mathcal
  O}_K[G]$-module (e.g. by ${\mathcal O}_K$-generators of ${\rm
  Sp}_J({G},{\mathcal O}_K)^I$, as was already used in the proof of Corollary
\ref{girr}), therefore there exists some $m>>0$ with $p_K^m{\rm
  Sp}_J({G},{\mathcal O}_K)\subset {\mathcal S}$. This means that (\ref{infini})
simplifies: it becomes ${\rm Sp}_J({G},{\mathcal O}_K)\subset p_K^n{\mathcal S}$. In view of this dichotomy (a)/(b) for any $n\in\mathbb{Z}$ we get $p_K^n{\mathcal S}={\rm Sp}_J({G},{\mathcal O}_K) $ for some $n\in{\mathbb Z}$ since $\bigcap_np_K^n{\mathcal S}=0$ and $\bigcup_np_K^n{\mathcal S}={\rm Sp}_J({G},K)$.\hfill$\Box$

\end{document}